\documentclass[a4paper,10pt]{article}

\usepackage[utf8]{inputenc}
\usepackage[T1]{fontenc}

\usepackage{a4wide}

\usepackage{amsmath,amssymb}
\usepackage{amsthm}
\usepackage{enumerate} 

\usepackage{authblk}

\usepackage[dvipsnames,table,xcdraw,svgnames]{xcolor}
\definecolor{highlightNEW}{named}{black}
\definecolor{myfullorange}{RGB}{255,128,0}
\definecolor{myfullred}{RGB}{244,63,43}
\definecolor{mygreen}{RGB}{70,180,5}
\definecolor{myhalforange}{RGB}{147,73,0}
\definecolor{myhalfred}{RGB}{158,22,7}
\definecolor{mylilas}{RGB}{170,55,241}
\definecolor{myred}{RGB}{244,63,43}

\usepackage{hyperref}
\hypersetup{
    breaklinks=true,
    colorlinks=true,
    linkcolor=Navy,
    citecolor=Navy,
    urlcolor=Navy
}

\usepackage{graphicx}
\graphicspath{{fig/}}
\usepackage{epstopdf,overpic}
\usepackage[section]{placeins} 

\usepackage{booktabs}
\usepackage{multirow}

\usepackage[absolute,overlay,showboxes]{textpos}
\TPMargin{2mm}
\textblockrulecolour{Navy}

\usepackage{natbib}

\usepackage{listings}
\usepackage{algorithm,algorithmic}

\newtheorem{theorem}{Theorem}[section]

\newtheorem{example}[theorem]{Example}

\newcommand{\doi}[1]{DOI~\href{\detokenize{http://dx.doi.org/#1}}{\detokenize{#1}}}
\newcommand{\zblnumber}[1]{Zbl~\href{\detokenize{https://zbmath.org/?q=an:#1}}{\detokenize{#1}}}
\newcommand{\mrnumber}[1]{\href{\detokenize{https://www.ams.org/mathscinet-getitem?mr=#1}}{\detokenize{MR#1}}}

\newcommand{\orcid}[1]{\includegraphics[scale=0.618]{icon_orcid_16x16.png}~\url{http://orcid.org/#1}}

\renewcommand{\d}{\,\mathrm{d}}

\newcommand{\ccode}[2]{\par
        \vspace*{8pt}
        {{\leftskip18pt\rightskip\leftskip
        \noindent{\it #1}\/: #2\par}}\par}
\newcommand{\keywords}[1]{\ccode{Keywords}{#1}}
\newcommand{\email}[1]{\href{mailto:#1}{#1}}

\def\received#1{Received~#1\par}
\def\revised#1{Revised~#1\par}
\def\accepted#1{Accepted~#1\par}
\def\published#1{Published~#1\par}

\def\foliofont{\fontsize{8}{10}\selectfont}
\usepackage[mathscr]{euscript}
\DeclareSymbolFont{rsfs}{U}{rsfs}{m}{n}
\DeclareSymbolFontAlphabet{\mathscrsfs}{rsfs}

\newcommand{\jpTitle}{Numerical aspects of integration in semi-closed option pricing formulas for stochastic volatility jump diffusion models}
\newcommand{\jpAuthors}{J. Dan\v{e}k and J. Posp\'{\i}\v{s}il}
\newcommand{\jpKeywords}{variable precision arithmetic; numerical integration; adaptive quadrature; option pricing; stochastic volatility models}
\newcommand{\jpMSC}{65D30; 91G60; 65R10}
\newcommand{\jpJEL}{C63}
\newcommand{\jpDateReceived}{23 October 2017} 
\newcommand{\jpDateRevised}{6 February, 1 June, 4 September, 12 December 2018}
\newcommand{\jpDateAccepted}{24 December 2018}
\newcommand{\jpDatePublished}{15 May 2019}
\newcommand{\jpDate}{}

\author[1]{Josef Dan\v{e}k} 
\author[1]{Jan Posp\'{\i}\v{s}il\thanks{Corresponding author, \email{honik@kma.zcu.cz}}} 
\affil[1]{NTIS - New Technologies for the Information Society, Faculty of Applied Sciences, \authorcr University of West Bohemia, Univerzitn\'{\i} 2732/8, 301 00 Plze\v{n}, Czech Republic,\vspace*{3pt}}

\ifpdf
\hypersetup{
  pdftitle={\jpTitle},
  pdfauthor={\jpAuthors},
  pdfkeywords={\jpKeywords},
  pdfinfo={
      MSC={\jpMSC},
      JEL={\jpJEL}
  }
}
\fi

\title{\textcolor{Navy}{\textsc{\jpTitle}}}
\date{\jpDate}

\begin{document}

\maketitle
\begin{textblock}{6}(1.25,1.15)
{\foliofont\noindent This is an Accepted Manuscript of an article published by Taylor \& Francis in the
International Journal of Computer Mathematics 
97(6), 1268--1292, 2020, \doi{10.1080/00207160.2019.1614174}. \\
Available online \url{https://www.tandfonline.com/10.1080/00207160.2019.1614174}.
}
\end{textblock}

\begin{center}
\received{\jpDateReceived}
\revised{\jpDateRevised}
\accepted{\jpDateAccepted}
\published{\jpDatePublished}
\end{center}

\begin{abstract}
In mathematical finance, a process of calibrating stochastic volatility (SV) option pricing models to real market data involves a numerical calculation of integrals that depend on several model parameters. This optimization task consists of large number of integral evaluations with high precision and low computational time requirements. However, for some model parameters, many numerical quadrature algorithms fail to meet these requirements. We can observe an enormous increase in function evaluations, serious precision problems and a significant increase of computational time.

In this paper we numerically analyse these problems and show that they are especially caused by inaccurately evaluated integrands. We propose a fast regime switching algorithm that tells if it is sufficient to evaluate the integrand in standard double arithmetic or if a higher precision arithmetic has to be used. We compare and recommend numerical quadratures for typical SV models and different parameter values, especially for problematic cases.

\end{abstract}

\ccode{MSC classification}{\jpMSC}
\ccode{JEL classification}{\jpJEL}

\clearpage

\section{Introduction}

Recently, \citet*{BaustianMrazekPospisilSobotka17asmb} presented a unifying approach to several stochastic volatility jump diffusion (SVJD) models. This approach among others covers the widely used \citet{Heston93} and \citet{Bates96} models, \citet{BarndorffNielsen01} model as well as a newly proposed approximative fractional stochastic volatility jump diffusion (AFSVJD) model \citep{PospisilSobotka16amf,MrazekPospisilSobotka16ejor}. Although we present the numerical pitfalls of numerical integration only for the AFSVJD model in detail here, similar numerical misbehaviour can be observed in all above mentioned models and probably in other SVJD models as well.

In many mathematical models we can observe that standard IEEE 32-bit or 64-bit floating-point arithmetic is not always fully sufficient, see for example works by \citet{Bailey05g}, who studied especially applications in physics \citep{Bailey12,Bailey15}. Computations that require higher than double precision for robust and exact decision making were introduced in \cite{Pal04}. To overcome the problems caused by the floating-point arithmetic limits, high-precision or variable-precision arithmetic is a rapidly growing part of scientific computing environments. In the above mentioned papers we can find among others currently available software packages for high-precision floating-point arithmetic. For example in MATLAB, there exists a possibility of defining variables and perform numerical calculations in variable-precision arithmetic using the \texttt{vpa} bundle that is part of the Symbolic Math Toolbox. Variable in the name suggests that a user can set the number of significant digits arbitrary, by default it is 32 significant digits. In this paper we show that in some cases it is actually necessary to use the \texttt{vpa} in order to get correct integration results and consequently correct option prices. 

Numerical integration in standard floating-point arithmetic is introduced in many university textbooks or monographs, let us mention at least the following books by \citet{Krylov62,Stoer02,Davis07} and by \citet{Dahlquist08}. The most common quadratures used in applications are the Gauss quadratures or the Simpson rule together with the adaptive refinement techniques. The problem of high-precision numerical integration is reviewed by \citet{Bailey11}. Gauss-Legendre quadratures using \texttt{vpa} in MATLAB was in particular studied by \citet{Rathod11}. In our case, we would like to get an effective, i.e. fast and sufficiently accurate, calculation of definite integral arising in the option pricing formula. This effectiveness can be for example achieved by clever switching algorithm between standard floating-point arithmetic that can be often sufficient and between \texttt{vpa}.

Although integration in option pricing models was already studied in several papers, according to authors' knowledge none of them focused on problems caused by inaccurately evaluated integrands or on high-precision integration. A~very good unpublished review of option pricing formulas based on Fourier transform is the online document by \citet{Schmelzle10}. Integrand in the Heston model was studied by \citet{KahlJackel05}. A~variations of the Fourier transform in option pricing was studied in \cite{Levendorskii12,BoyarchenkoLevendorskii14}, however, from the numerical point of view the usage of trapezoidal rule is far from being satisfactory. A~widely used techniques in option pricing are based on Fourier method. Among these methods we find the classical fast Fourier transform (FFT) as was suggested by \citet{CarrMadan99} or the fractional FFT modification \citep{Bailey91,Bailey94}, Fourier method with the Gauss-Laguerre quadrature \citep{Lindstrom08}, the so called COS method \citep{Fang08} or methods based on wavelets \citep{Ortiz16}. All these methods can be fast for example in calculating an approximation of the integral in many discrete points at once, however, in option pricing problems many values are calculated redundantly and moreover with relatively low precision that should be in modern financial applications considered unsatisfactory.

The paper is structured as follows. 
In Section \ref{sec:svmodels} we introduce the studied formula for European call option price obtained by the approximative fractional stochastic volatility model. 

The problem of inaccurately evaluated integrand is presented in Section \ref{sec:integrand}. A~special attention is paid to study the integral behaviour during the optimization process that occur during the calibration of the model to real market data. We show that an inaccurately evaluated integrand can lead to changes in option price in order of hundreds of dollars. We suggest a usage of the variable precision arithmetic for problematic cases and design a switching regime algorithm.

In Section \ref{sec:quadratures} we explain why numerical quadratures fail.  
Further we compare several numerical quadratures in Section \ref{sec:results} where we give also some recommendations what quadrature to use and why or why not, how to handle adaptivity and 
how to set the tolerances to get reliable results. In particular we compare the calibration results for the cases where the switching regime algorithm is used or not.
We conclude in Section~\ref{sec:conclusion}.

\section{Stochastic volatility models}\label{sec:svmodels}

Following \citet*{BaustianMrazekPospisilSobotka17asmb}, we consider a general stochastic volatility jump diffusion (SVJD) model that covers several kinds of stochastic volatility processes and also different types of jumps
\begin{align*}
\d S_t &= (r-\lambda\beta)S_t \d t + \sqrt{v_t}S_t \d W^S_t + S_{t-} \d Q_t,\\
\d v_t &= p(v_t) \d t + q(v_t) \d W^v_t,\\
\d W^S_t \d W^v_t &= \rho\d t,
\end{align*}
where $p,q\in C^{\infty}(0,\infty)$ are general coefficient functions (for particular choices of $p$ and $q$ see Table \ref{t:pq}),
$r$ is the interest rate, $\rho$ is the correlation of Wiener processes $W_t^S$ and $W_t^v$,
parameters $\lambda$ and $\beta$ correspond to a specific jump process $Q_t$, which is
a compound Poisson process $Q_t = \sum\limits_{i=1}^{N_t} Y_i$, where
$Y_1, Y_2, \dots$ are pairwise independent random variables with identically distributed jump sizes $\beta=\mathbb{E}[Y_i]$ for all $i\in\mathbb{N}$, $N_t$ is a standard Poisson process with intensity $\lambda$ independent of the $Y_i$.

\begin{table}[ht!]
\caption{Different SVJD models}\label{t:pq}
\begin{center}
\begin{tabular}{lll}
model & $p(v)$ & $q(v)$ \\
\toprule
Heston/Bates & $\kappa(\theta-v)$ & $\sigma\sqrt{v}$ \\
3/2 model$^\ast$ & $\omega v - \tilde{\theta}v^2$ & $\xi v^{\frac{3}{2}}$ \\
Geometric BM & $\alpha v$ & $\xi v$ \\
AFSVJD$^{\ast\ast}$ & $(H-1/2)\psi_t\sigma\sqrt{v}+\kappa(\theta-v)$ & $\varepsilon^{H-1/2}\sigma\sqrt{v}$ \\
\bottomrule
\end{tabular}\\
\begin{flushleft}
\hspace*{42mm}
$^\ast${\footnotesize $\tilde{\theta}=-\frac{1}{2}\xi^2+(1-\gamma)\rho\xi+\sqrt{(\theta+\frac{1}{2}\xi^2)^2-\gamma(1-\gamma)\xi^2}$,} \\
\hspace*{42mm}
$^{\ast\ast}${\footnotesize $\psi_t= \int_0^t (t-s+\varepsilon)^{H-3/2} \d W^{\psi}_s$. }
\end{flushleft}
\end{center}
\end{table}

A \emph{unifying} formula for the price $V=V(K,\tau)$ of a European call option with strike price $K$ and time to maturity $\tau$ has the form \citep{BaustianMrazekPospisilSobotka17asmb}
\begin{equation}\label{e:unif}
V = S - Ke^{-r\tau}\frac1{2\pi} \int_{-\infty+ik_i}^{\infty+ik_i} e^{-ik\tilde{X}} e^{\lambda(\hat{\varphi}(-k)-1)\tau} \frac{\hat{F}(k,v,\tau)}{k^2-ik} \d k, 
\end{equation}
where $\tilde{X}= \ln(S/K) + (r-\lambda\beta)\tau$ and $\max(k_1,0)<k_i<\min(1,k_2)$ and
$\hat{F}$ is the so called \emph{fundamental transform} of the particular stochastic volatility part and
$\hat{\varphi}$ is the transform of the jump term. An integration domain in the complex plane is a line (represented by $k_i$) lying in the suitable \emph{strip of regularity} \cite{Lewis00,Lewis16}, for European call options it suffices to take $k_i = 1/2$.

If $B^{\varepsilon}_t = \int\limits_0^t (t-s + \varepsilon)^{H-\frac12} \d W_s$ is the \emph{approximative fractional Brownian motion},
$\varepsilon>0$ (for $\varepsilon\to 0$ the process converges to the \emph{standard fractional Brownian motion}), 
$H>1/2$ (for $H=1/2$ it is the \emph{standard Brownian motion}), 
then the volatility process in the \emph{approximative fractional SVJD} (AFSVJD) model \cite{PospisilSobotka16amf}
\begin{align*}
\d v_t &= \kappa (\theta - v_t) \d t + \sigma \sqrt{v_t} {\d B^{\varepsilon}_t}, 
\intertext{can be rewritten as}
\d v_t &= \left[(H-1/2)\psi_t\sigma\sqrt{v_t} + \kappa (\theta - v_t)\right] \d t + \varepsilon^{H-1/2}\sigma \sqrt{v_t} \d {W}^v_t,
\end{align*}
where $\psi_t = \int_0^t (t-s+\varepsilon)^{H-3/2} dW^{\psi}_s$.

\begin{example}
Jumps examples:
\begin{enumerate}
\item In \citet*{Bates96} model, jump sizes are log-normal, 
$\ln(1+Y_i) \sim \mathcal{N}(\mu_J,\sigma_J^2)$, \\
$\hat{\varphi}(k) = \exp\left\{ i\mu_J k -\frac12 \sigma_J^2 k^2 \right\}$,
$\beta = \hat{\varphi}(-i)-1 = \exp\left\{ \mu_J + \frac12 \sigma_J^2 \right\}-1$. \\[3pt]
\item In \citet*{YanHanson06} model, jump sizes are log-uniform,
$\ln(1+Y_i) \sim \mathcal{U}(a,b)$, \\
$\hat{\varphi}(k) = \frac{e^{ikb}-e^{ika}}{(b-a)ik}$,
$\beta = \hat{\varphi}(-i)-1 = \frac{e^b-e^a}{b-a}-1$.
\end{enumerate}
\end{example}

In order to perform a thorough numerical analysis of the integral \eqref{e:unif} a particular model with particular fundamental transform $\hat{F}$ has to be considered. Sometimes it is also useful to represent the jumps in terms of a \emph{characteristic function} that we denote by $\phi$. Price $V=V(K,\tau)$ of a European call option with strike price $K$ and time to maturity $\tau$ in the AFSVJD model is given \citep{MrazekPospisilSobotka16ejor,BaustianMrazekPospisilSobotka17asmb} by 

\begin{equation}
V = S - Ke^{-r\tau}\frac{1}{\pi}\int\limits_{0+i/2}^{+\infty+i/2} \underbrace{e^{-ikX}\frac{\hat{F}(k,v,\tau)}{k^2-ik}\phi(-k)}_{f(k)} dk,\label{e:price}
\end{equation}
where $X = \ln({S}/{K}) + r\tau$, fundamental transform 
\begin{align}
\hat{F}(k,v,\tau) &= \exp (C(k,\tau) + D(k,\tau) v),\label{e:hatF} \\
C(k,\tau) &= \kappa\theta \Biggl( Y\tau - \underbrace{\frac{2}{B^2}}_{C_1} \underbrace{\ln\left(\frac{1-ge^{-d\tau}}{1-g}\right)}_{C_2} \Biggr),\notag\\
D(k,\tau) &= Y \frac{1-e^{-d\tau}}{1-ge^{-d\tau}},\notag\\
Y &= -\frac{k^2-ik}{b+d},\quad g = \frac{b-d}{b+d},\notag\\
d &= \sqrt{b^2+B^2(k^2-ik)},\notag\\
b &= \kappa + ik\rho B,\quad B = \varepsilon^{H-1/2}\sigma,\notag
\intertext{and characteristic function}
\phi(k) &= \exp\left\{ -i \lambda \beta k~\tau + \lambda \tau \biggl[\hat{\varphi}(k)-1 \biggr] \right\}, \label{e:phi}\\
\beta &= \exp\left\{ \mu_J + \frac12 \sigma_J^2 \right\}-1,\notag\\
\hat{\varphi}(k) &= \exp\left\{ i\mu_J k~-\frac12 \sigma_J^2 k^2 \right\}.\notag
\end{align}

In Table \ref{tab:bounds} we can see typical simple lower and upper real-valued bounds (LB and UB) for model parameters 
$\chi = (v_0, \kappa, \theta, \sigma, \rho, \lambda, \mu_J, \sigma_J, H)$
considered in calibration to real market data. It is worth to mention that for $H=0.5$ the AFSVJD model coincides with the \citet{Bates96} model and if further there are no jumps ($\lambda=0$), we get the \citet{Heston93} model. For this reason, AFSVJD model was chosen as a particular model to study. In \cite{BaustianMrazekPospisilSobotka17asmb}, authors also showed that a similar type of integral can be obtained also for models that do not follow exactly the above mentioned stochastic dynamics such as the \citet*{BarndorffNielsen01} model. Similar pricing formulas for a variety of other SVJD models such as models with other L\'{e}vy processes (double exponential, variance gamma, normal inverse Gaussian, normal tempered stable, finite moment log stable) or with different type of stochastic volatility (SABR, 3/2 model) can be found in the book by \citet*{Lewis16}. 

\begin{table}[h]
\caption{Typical simple lower and upper real-valued bounds (LB and UB) for model parameters considered in calibration process.}
\label{tab:bounds}
\centering
\begin{tabular}{lccccccccc}
   & $v_0$ & $\kappa$ & $\theta$ & $\sigma$ & $\rho$ & $\lambda$ & $\mu_J$ & $\sigma_J$ & $H$ \\
\midrule
LB: & 0     & 0        & 0        & 0        & -1     & 0         & -10     & 0          & 0.5 \\
UB: & 1     & 150      & 1        & 4        & 1     & 100       & 5       & 4          & 1.0 \\
\bottomrule
\end{tabular}
\end{table}

In our experiments we consider the following ranges that are typical for many option market data sets (options on major indices like DAX, FTSE 100, Nikkei 225 and S\&P 500 that occurred at the market in the last two years):%
\begin{description}
\item[price] of the underlying asset (usually in dollars): $0< S\leq S_{\max}=30\ 000$,
\item[time to maturity] (in years): $0< \tau\leq \tau_{\max}=5$, 
\item[strike price]:               $0< K\leq K_{\max}=3S_{\max}$,
\item[interest rate] (positive):   $0< r\leq r_{\max}=0.05$.
\end{description}
We will denote the vector of market data by $\psi=(\tau,K,r,S)$.

A~process of calibrating option pricing models to real market data involves a~numerical calculation of integrals similar to the integral in \eqref{e:price}. As we can see, the value of the integral depend highly nonlinearly on several model parameters and real market data. Calibration as an optimization task consists of large number of integral evaluations that must be calculated with high precision, but at low computational costs. However, for some model parameters, many numerical quadrature algorithms fail to meet these two requirements. We can observe an~enormous increase in function evaluations (especially in adaptive quadratures), serious precision problems (even with the simple non-adaptive trapezoidal rule) as well as a~significant increase of computational time. 

\section{Inaccurately evaluated integrand}\label{sec:integrand}

In this section we focus on the deeper numerical analysis of the above mentioned behaviour. We show that some of the problems are caused by an inaccurate evaluation of the integrand in standard \texttt{double} precision arithmetic \citep{DanekPospisil15tcp}. In the IEEE Standard for Floating-Point Arithmetic (IEEE 754), the smallest interchange format for the standard \texttt{double} precision number is referred as \texttt{binary64} and it contains 1 sign bit, 11 exponent bits and 53 significand (or mantissa) precision bits (only 52 are explicitly stored). In \texttt{double} precision arithmetic we therefore get 15--17 significant decimal digits precision. Moreover, even simple arithmetic operations such as addition and subtraction can lead to loosing significant decimal digits precision.

\begin{example}[Loosing significant decimal digits precision]\label{ex:loose_digits}
Let us consider a simple example $10^{-1} + 10^9 -10^9$ whose result in \texttt{double} arithmetic is 0.100000023841858, i.e. after the two operations the number of significant digits lowered to half. The problem is already with the number 0.1 that has periodic representation in binary: $0.1 = (0.0\overline{0011})_2$, i.e. 0.1 in \texttt{double} is
\[ 0.\underbrace{1100\,1100 \dots 1100}_{48\text{ digits}} 11010 \times 2^{-3}. \]
On the other hand, the number $10^9$ has an exact representation in \texttt{double}
\[ 0. 111011100110101100101 \times 2^{30}. \]
Addition of these two numbers shifts the exponent of the lower number to the higher exponent, i.e. 0.1 after the shift
\[ 0.\underbrace{0000 \dots 0}_{33\text{ digits}} 1100\,1100\,1100\,1100\,1101 \times 2^{30} \]
which is in fact the resulting number after adding and subtracting $10^9$, in decimal 0.100000023841858.
\end{example}

Loosing significant digits consequently leads to inaccurate evaluation of functions.

\begin{example}[Inaccurate evaluation of a simple function]\label{ex:inexact_simple}
Let $h(x)=\delta x^2$, where $\delta$ is a given real parameter, say $10^6$. Its implementation in MATLAB can be for example

\begin{lstlisting}
 function result = h_exact(x);
   delta  = 1e6; 
   result = delta*x.^2;
 end
\end{lstlisting}

The following function should theoretically return the same values. 

\begin{lstlisting}
 function result = h_trouble(x);
   delta  = 1e6; 
   a      = x + 1000*delta;
   b      = (a.^2+delta*(x.^2)).^0.5;
   result = b.^2-a.^2;
 end
\end{lstlisting}

However, in the view of the previous Example \ref{ex:loose_digits}, we can expect problems with loosing significant digits precision. Indeed, in Figure \ref{fig:inexact_simple} there are graphs of both functions depicted over the interval $[0,1]$ or in a neighbourhood of number $0.5$. Function \texttt{h\_exact(x)} is \textbf{\textcolor{blue}{bold blue}} and function \texttt{h\_trouble(x)} is \textcolor{red}{red}. Although the global view of the function $h(x)$ does not indicate any troubles, a~more detailed view (zoom) reveals unexpected discontinuities in the function that should be smooth. 

\begin{figure}[!ht]
\centering
\includegraphics[width=\textwidth]{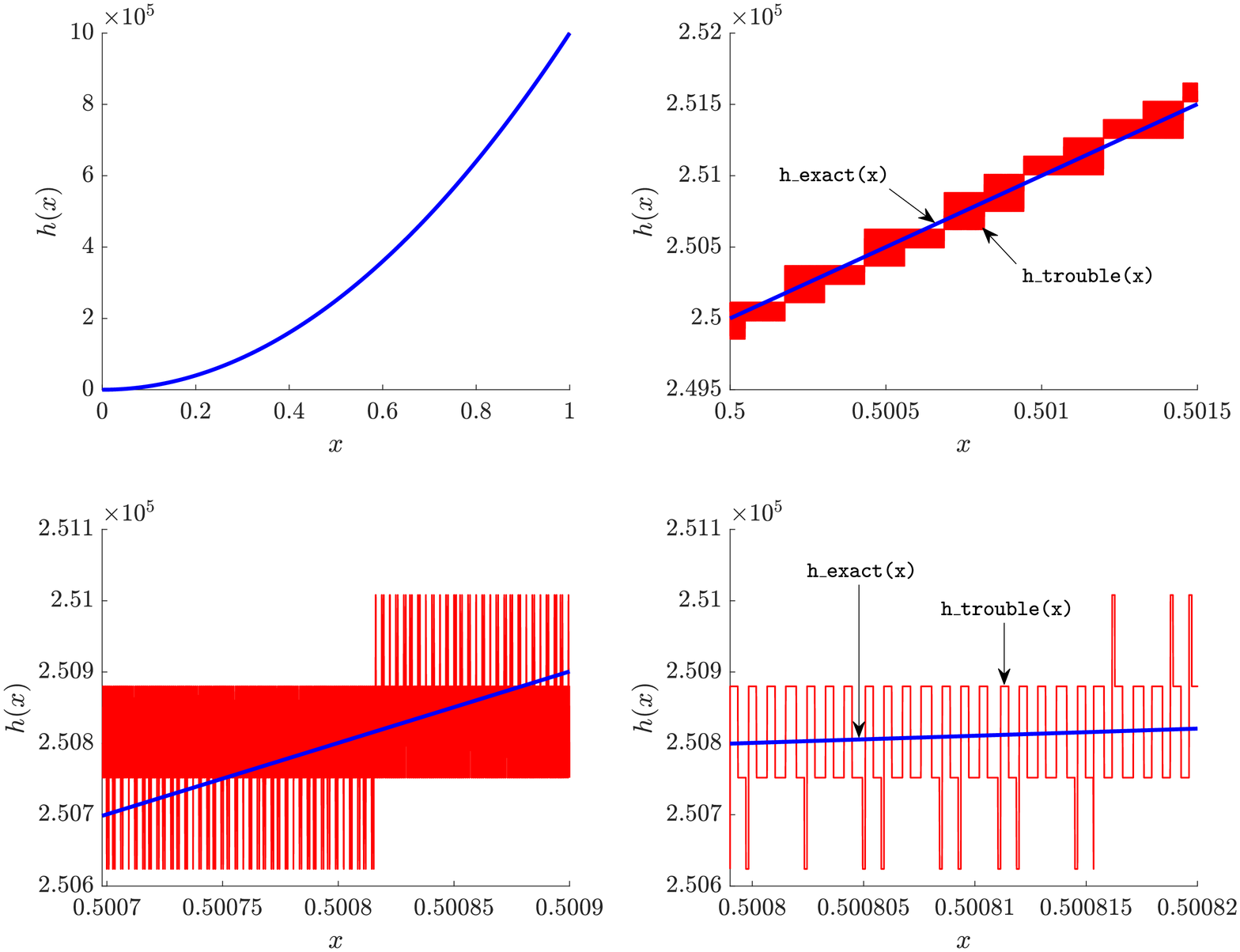}
\caption{Global and detailed view to \texttt{h\_exact(x)} and \texttt{h\_trouble(x)} from Example \ref{ex:inexact_simple}.}\label{fig:inexact_simple}
\end{figure}

\end{example}

In Example \ref{ex:inexact_simple}, it is clear that the form of \texttt{h\_trouble(x)} was implemented inefficiently and one should use the simpler form \texttt{h\_exact(x)}. However, in more complicated examples such a simplification is not always possible. The function of interest is the integrand $f(k)$ in \eqref{e:price}. In the following we show that if it is inaccurately evaluated, we can observe similar misbehaviour as in the previous simple example. To avoid the problems with loosing significant digits, we perform the evaluation of function values also in the high precision arithmetic, in particular in \texttt{vpa} in MATLAB. All \texttt{vpa} values in this paper are obtained with 32 significant decimal digits precision. 

From now on, the numerical analysis will involve only the AFSVJD model and the formula \eqref{e:price}. In other SVJD models, very similar analysis can be performed with analogous switching regime setting (see below).

\begin{figure}[ht!]
\centering
\begin{overpic}[width=\linewidth]{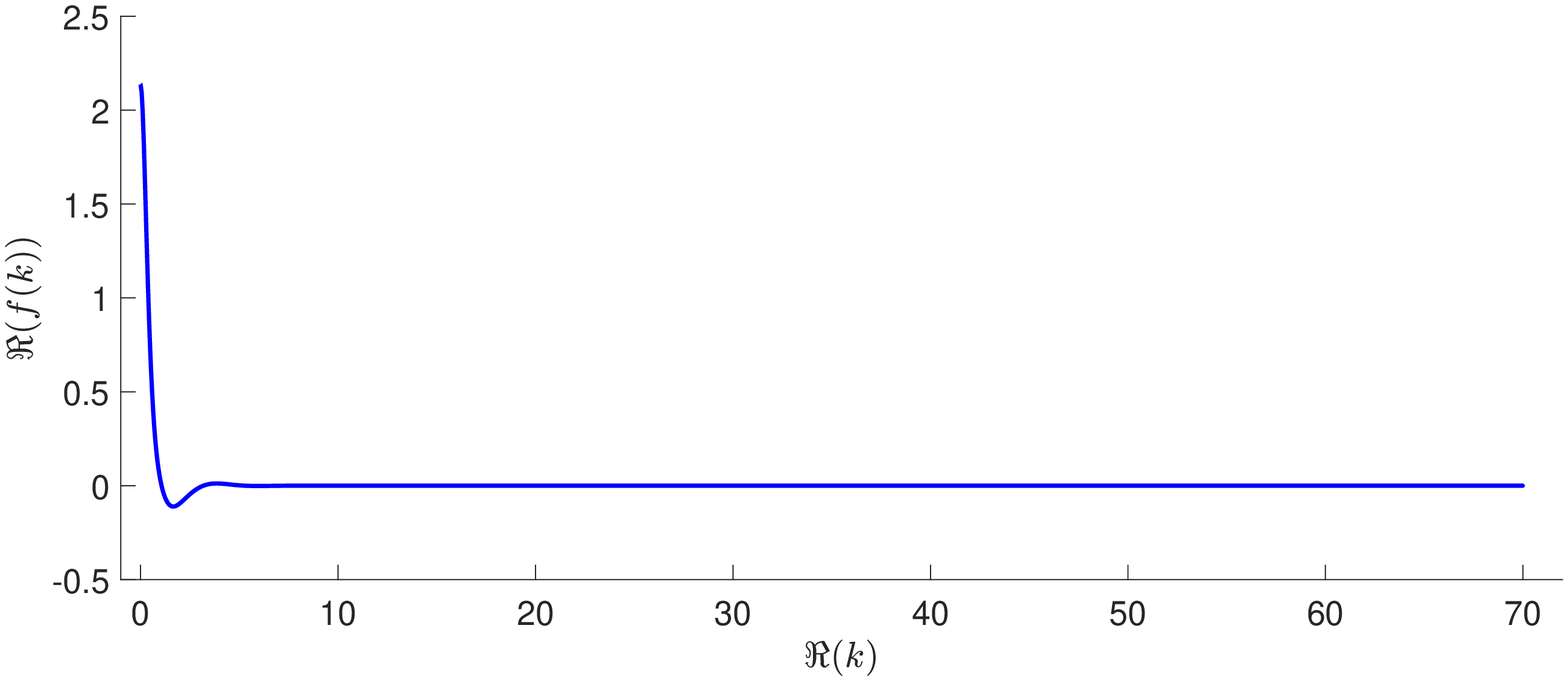}
\put(25,40){\small $v=0.97$, $\kappa=17.6$, $\theta=0.95$, $\rho=-0.86$}
\put(25,35){\small $\lambda=11.7$, $\mu_J=-6.66$, $\sigma_J=1.007$, $H=0.96$}
\end{overpic}
\\[5mm]
\begin{overpic}[width=0.5\linewidth]{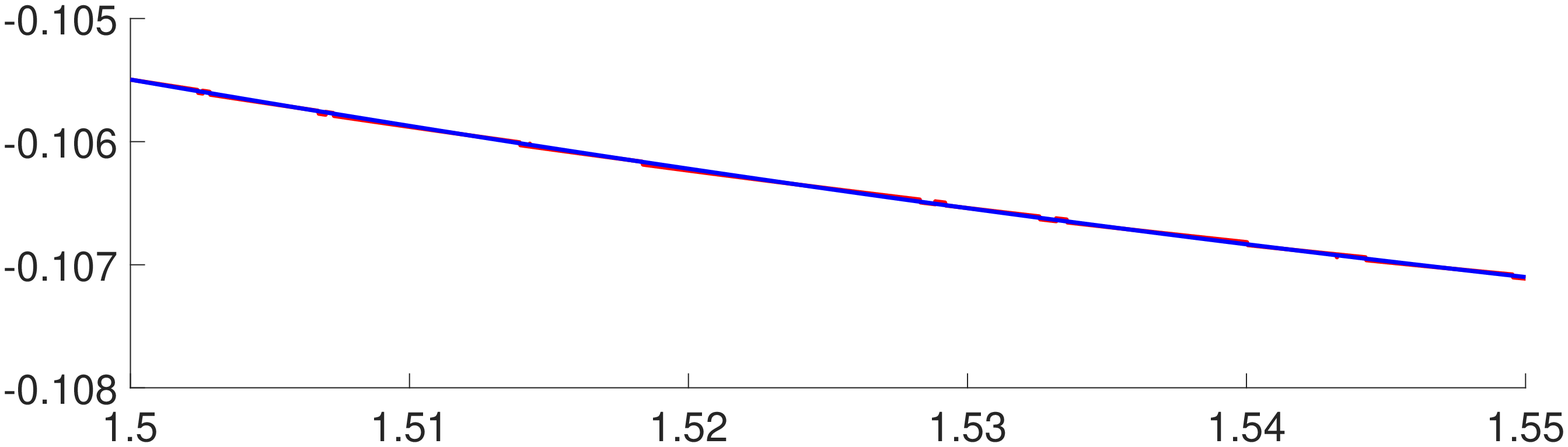}\put(40,24){\small $\sigma=0.0001$}\end{overpic}%
\begin{overpic}[width=0.5\linewidth]{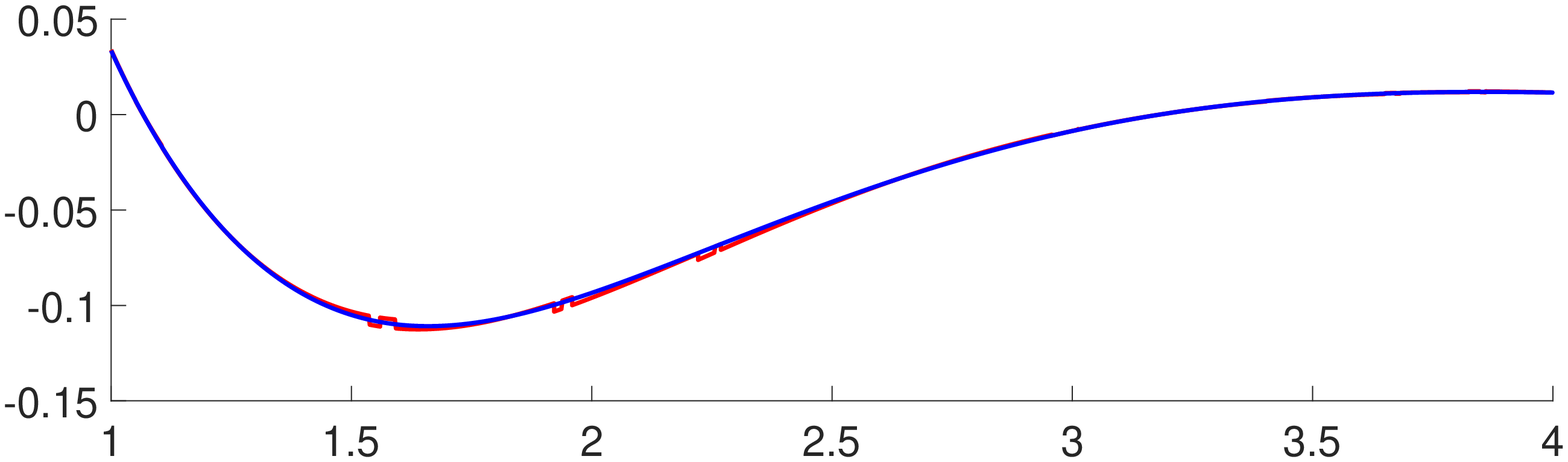}\put(40,24){\small $\sigma=0.00001$}\end{overpic}\\[3mm]
\begin{overpic}[width=0.5\linewidth]{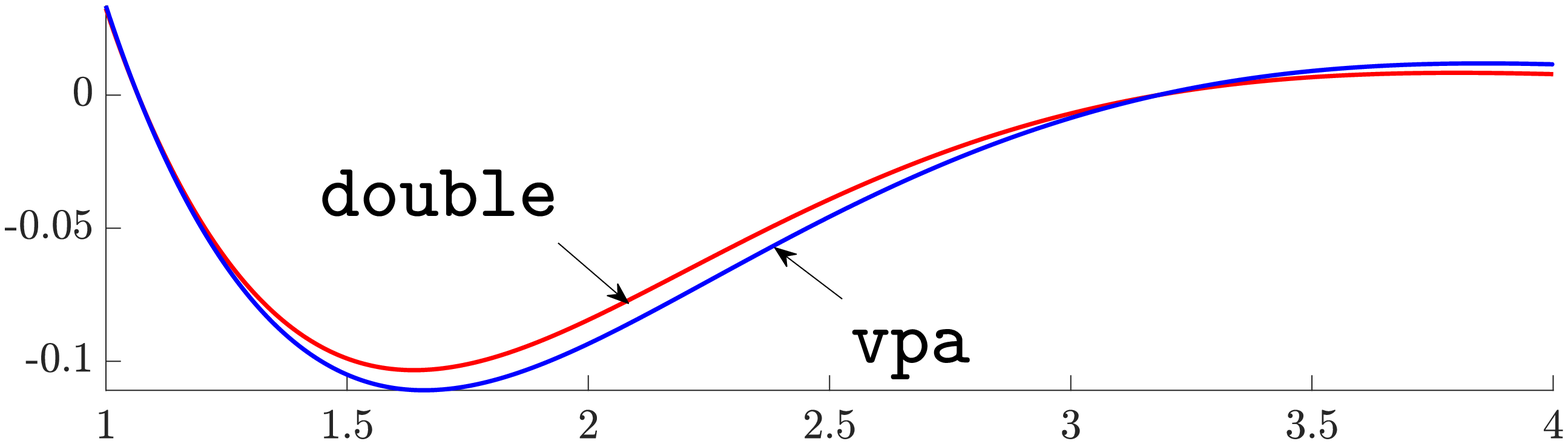}\put(40,24){\small $\sigma=0.000001$}\end{overpic}%
\begin{overpic}[width=0.5\linewidth]{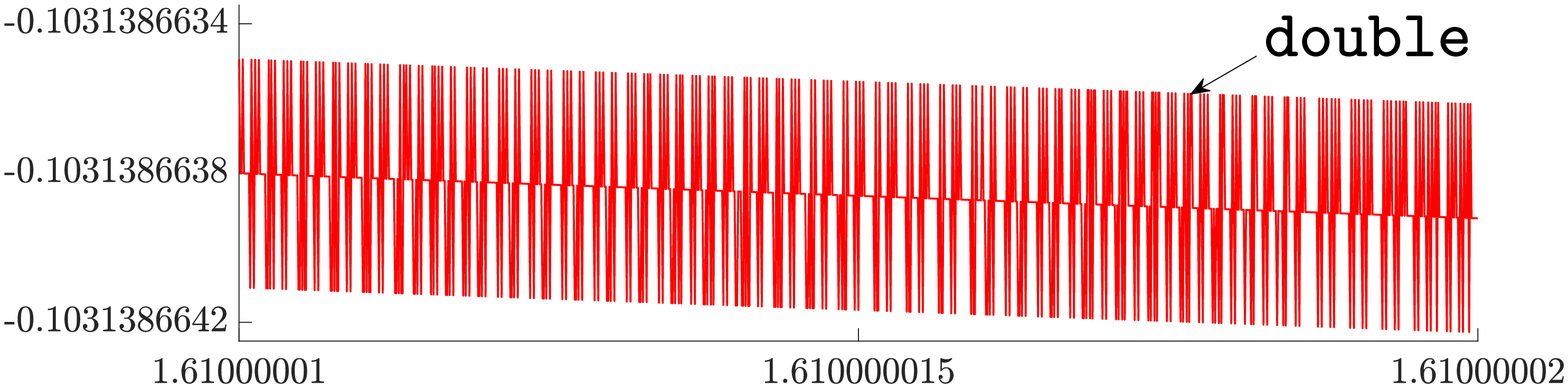}\put(40,24){\small $\sigma=0.000001$}\end{overpic}
\caption{Global and local view to the integrand $f(k)$ in Test Case 1.}\label{fig:test1}
\end{figure}

\begin{figure}[ht!]
\centering
\begin{overpic}[width=\linewidth]{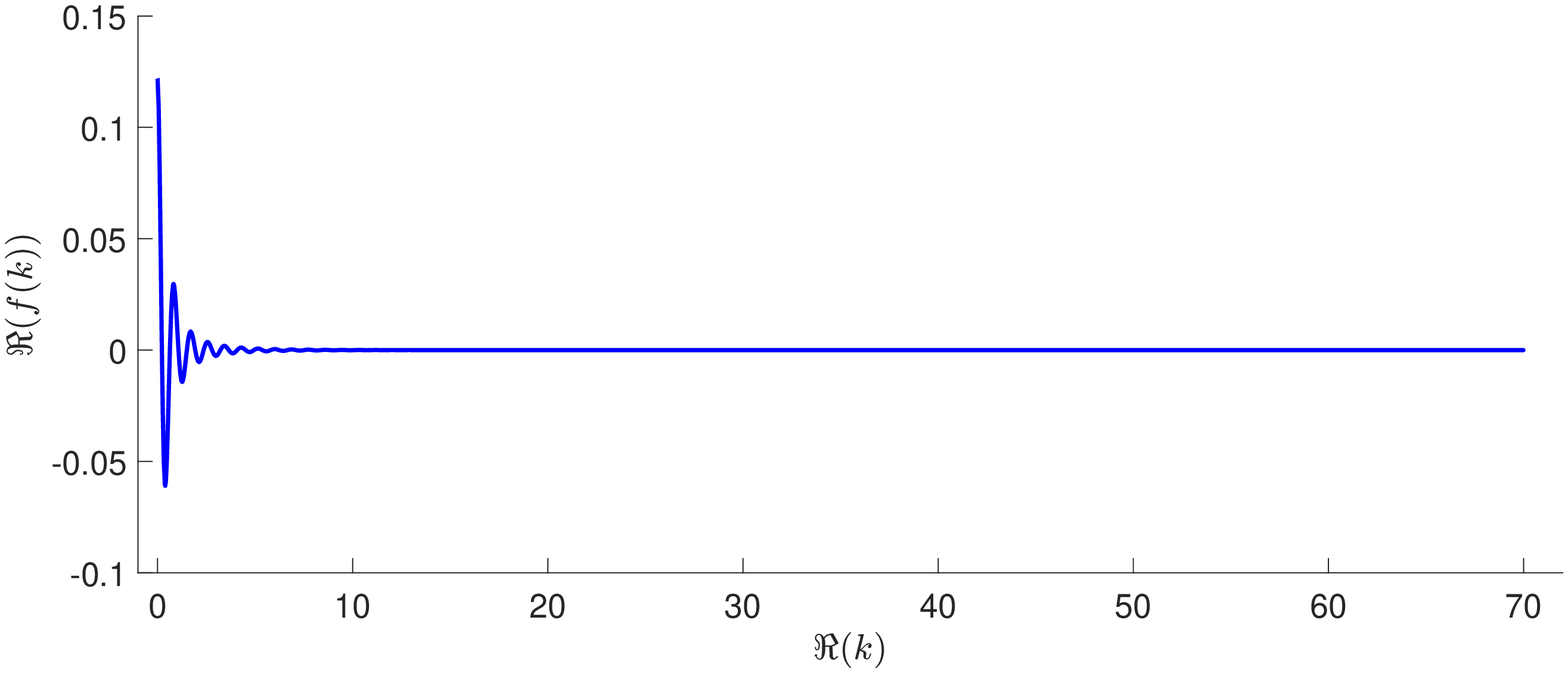}
\put(25,40){\small $v=0.3$, $\kappa=5$, $\theta=0.1$, $\rho=-0.5$}
\put(25,35){\small $\lambda=60$, $\mu_J=-9$, $\sigma_J=1.1$, $H=0.6$}
\end{overpic}
\\[5mm]
\begin{overpic}[width=0.5\linewidth]{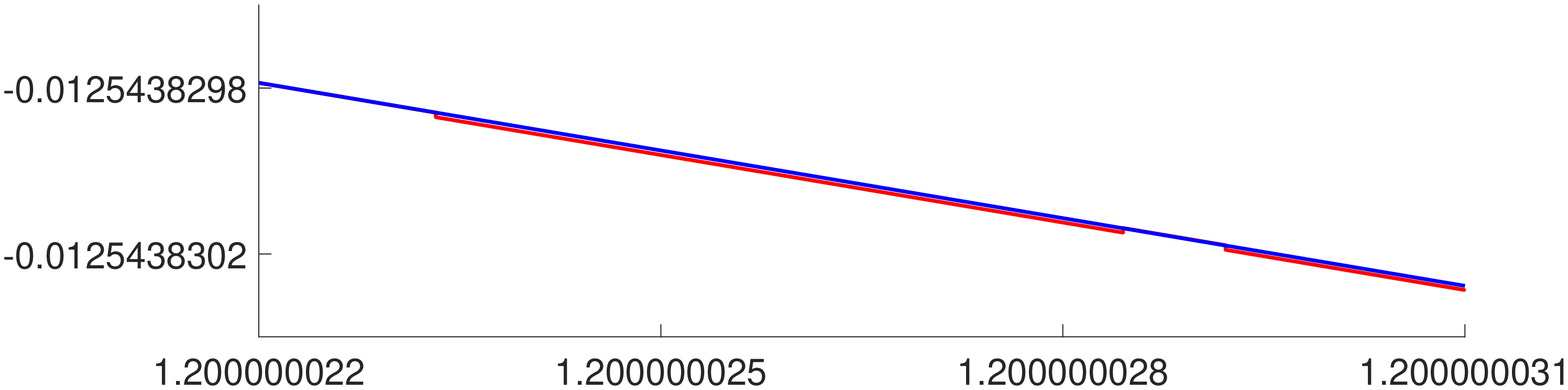}\put(40,24){\small $\sigma=0.001$}\end{overpic}%
\begin{overpic}[width=0.5\linewidth]{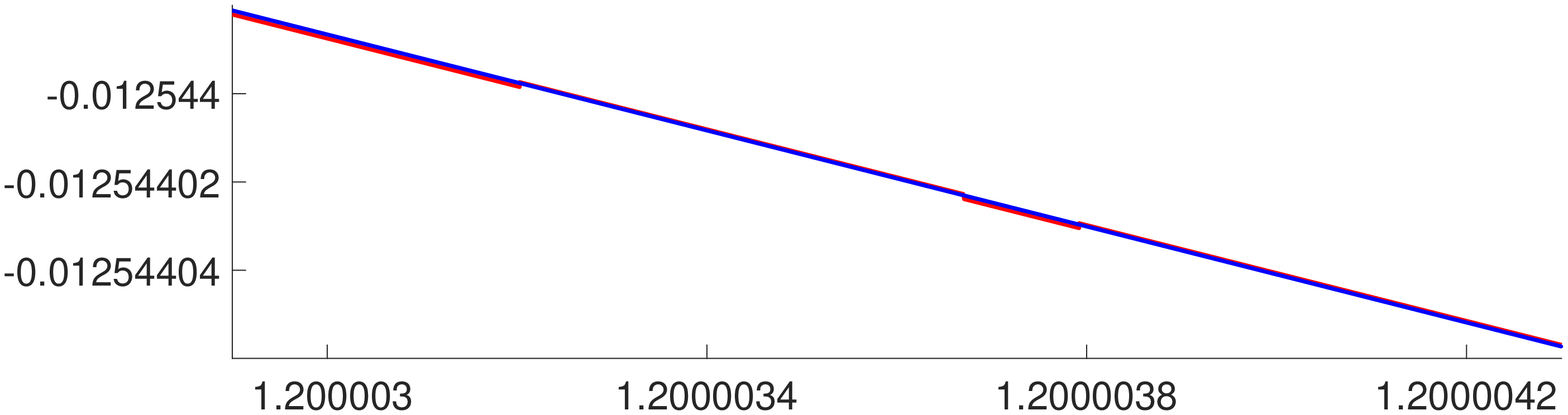}\put(40,24){\small $\sigma=0.0001$}\end{overpic}\\[3mm]
\begin{overpic}[width=0.5\linewidth]{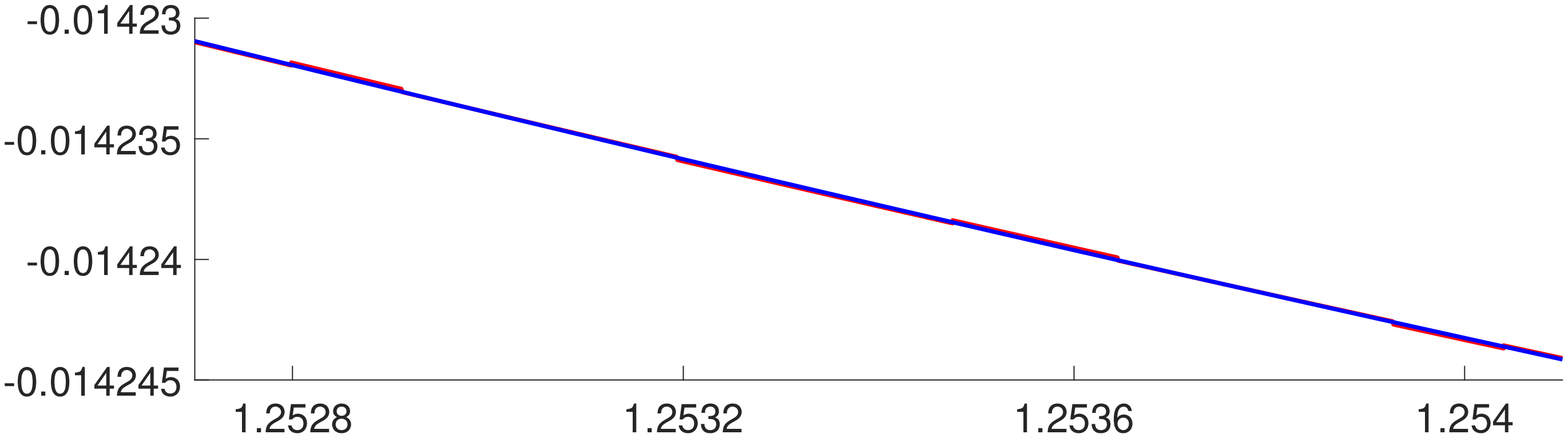}\put(40,24){\small $\sigma=0.00001$}\end{overpic}%
\begin{overpic}[width=0.5\linewidth]{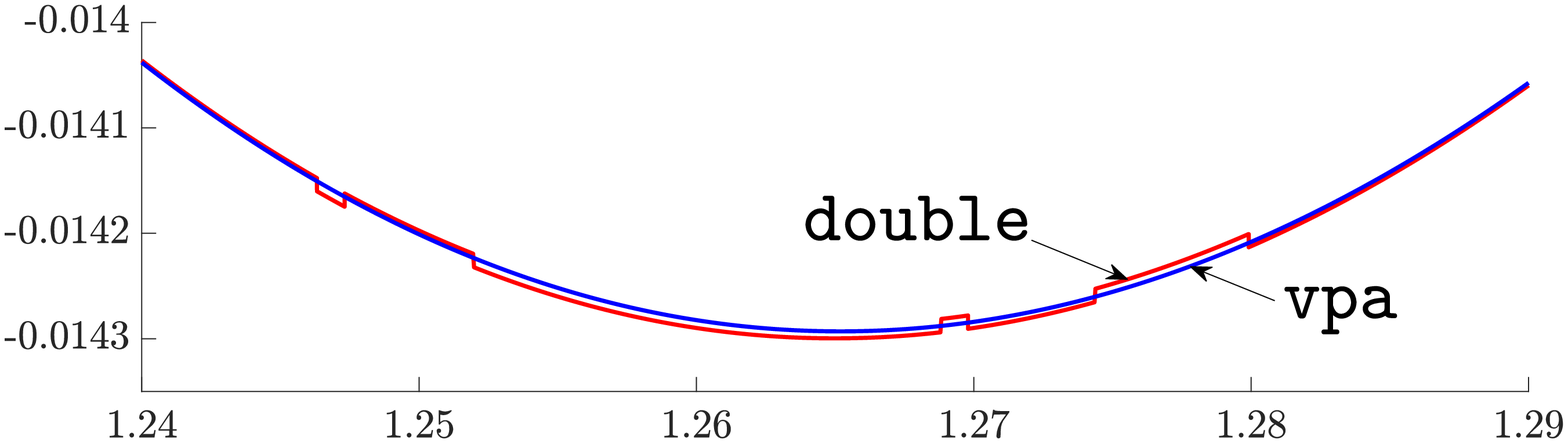}\put(40,24){\small $\sigma=0.000001$}\end{overpic}
\caption{Global and local view to the integrand $f(k)$ in Test Case 2.}\label{fig:test2}
\end{figure}

\begin{example}[Global and local view to the integrand $f(k)$]\label{ex:test_cases}
In Figures \ref{fig:test1} and \ref{fig:test2} we consider market data (European call options to FTSE 100 dated 8 January 2014, see \citep{PospisilSobotka16amf} and Example \ref{ex:calib} below) 
$$\psi = (0.120548, 6250, 0.009, 6721.8)$$ and demonstrate the numerical misbehaviour by changing only values of $\sigma$ and also the zoom depth, other parameters remain in each example the same.
Whereas the inaccurately enumerated values are in \texttt{double} precision arithmetic (\textcolor{red}{red} in Figures \ref{fig:test1} and \ref{fig:test2}), the smooth values are in \texttt{vpa} (\textcolor{blue}{blue}). In Figure \ref{fig:test1}, bottom pictures for $\sigma=0.000001$, we can see that potential smoothing of the inaccurately evaluated integrand would not help, since it would be completely set off the exact values. Bottom right picture hence shows only the zoom of the \texttt{double} evaluated integrand from the bottom left picture.
\end{example}

Moreover, inaccurately evaluated integrand can lead to really big differences in the option price, even in the order of hundreds of dollars. 
\begin{example}[Hundred dollars error]\label{ex:hundred_dollars}
For parameters 
\begin{align*}
\chi &= (0.98, 8, 0.8, 0.000001, -0.75, 0.75, 1.4, 0.2, 0.9)
\intertext{that can be observed during a calibration process and for data}
\psi &= (0.34, 12 500, 0.017, 10 000)
\end{align*}
we get
\begin{verbatim}
price_double=4115.317   price_vpa=3999.167   error=116.150
  int_double=1.48755534   int_vpa=1.51691623 error=0.02936088
\end{verbatim}
where \texttt{price} denotes the value of the option price calculated using the Gauss-Kronrod(7,15) quadrature that is by default used in the function \texttt{integral} in MATLAB (we compare other numerical quadratures later in Section \ref{sec:quadratures}) and \texttt{int} denotes the value of the integral. The sub indices \verb|_double| and \verb|_vpa| denote that the value was obtained for integrand $f(k)$ evaluated using the standard \texttt{double} precision or \texttt{vpa} arithmetic respectively. 
Absolute difference \texttt{error} indicate a huge difference in the option price value.
\end{example}

In order to get a reasonable values in the calibration process, we need the integral to be calculated with the error less then \texttt{1e-8}. The reason for this tolerance is the fact that in one calibration procedure, several hundred options are typically considered, which means that the utility (or fitness) function in the optimization procedure consists of this amount of integral evaluations and hence the option price calculation in each iteration. A safe option price value tolerance can be considered \texttt{1e-3}. From the comparison of the error in integration and option price in the example above we can clearly see that in order to get the safe price within error \texttt{1e-3}, we need integral to be calculated with error less then \texttt{1e-8}.

It is worth to mention that the integrand $f(k)$ does not have to be evaluated in \texttt{vpa} for all parameter values. In fact, the most problematic cases occur only for rare certain combinations of parameters. In Table \ref{tab:2-1-cetnosti} we present an influence of one model parameter (in particular parameter $\sigma$) and approximation parameter $\varepsilon$ to the number of problematic cases. A problematic case is such that the integral was evaluated either with error greater than \texttt{1e-8}, or with \texttt{fevals} greater than \texttt{1e4} (i.e. with high CPU time). Among all problematic cases, more than a third gave an error greater than \texttt{1e-2}. Values in Table \ref{tab:2-1-cetnosti} show number of problematic cases in 10 millions random integral evaluations, parameter $\sigma$ is generated uniformly in different ranges with approximation parameter $\varepsilon$ being fixed in all corresponding trials. All other parameters are chosen uniformly random in bounds from Table \ref{tab:bounds}. Influence of all parameters to the problematic cases follows below.

\begin{table}[t!]
\caption{Number of problematic cases in 10 millions random integral evaluations for different ranges for parameter $\sigma$ and different values of approximation parameter $\varepsilon$. All other parameters are chosen uniformly random in bounds from Table \ref{tab:bounds}.}\label{tab:2-1-cetnosti}
\begin{center}
\begin{tabular}{lrrrr}
$\sigma\in\:\backslash\:\varepsilon =$ 
                    &   $10^{-3}$ &    $10^{-4}$ &    $10^{-5}$ &    $10^{-6}$ \\
\midrule
$[10^{-5};4]$       &          42 &          103 &          250 &          696 \\
\midrule
$[10^{-2};10^{-1}]$ &           0 &            0 &          155 &       7\,315 \\
$[10^{-3};10^{-2}]$ &         262 &      10\,742 &      65\,979 &     179\,818 \\
$[10^{-4};10^{-3}]$ &    109\,922 &     269\,021 &     418\,708 &     528\,910 \\
$[10^{-5};10^{-4}]$ &    684\,565 &     782\,732 &     841\,482 &     882\,383 \\
$[10^{-6};10^{-5}]$ & 1\,059\,632 &  1\,064\,681 &  1\,068\,015 &  1\,071\,241 \\
\bottomrule
\end{tabular}
\end{center}
\end{table}

To overcome the raised issues, we introduce a switching regime for problematic cases, i.e. in majority of cases, the evaluation will be done in standard \texttt{double} arithmetic, however, in problematic cases \texttt{vpa} arithmetic has to be used. 

When analysing the integrand $f(k)$ in \eqref{e:price}, we can find out that the order of magnitude in terms $C_1$ and $C_2$ can differ significantly. In fact, this difference is the core of the problem of inaccurately evaluated integrand. It is worth to mention that presented form of $C$ is the most numerically stable. Although it is possible to put the term in front of the logarithm into the exponent of the argument, which is mathematically equivalent, a presented form is numerically more suitable for finite precision arithmetic calculations. We introduce the following algorithm for the switching between the \texttt{double} and \texttt{vpa} evaluation of the integrand.

Let $\Re(z)$ denote the real part of a complex number $z$. Let $o_1$ and $o_2$ denote the order for the terms $C_1$ and $C_2$ respectively and $f_0$ be the value 
\begin{equation}\label{e:f0}
f_0 := |\Re (f(0+i/2))|.
\end{equation}
Then 
\begin{align*}
o_1 &:= \log_{10}|\Re (C_1)|, \\
o_2 &:= \log_{10}|\Re (C_2)|, \\
o   &:= o_1 - o_2,
\end{align*}
where we consider $\log_{10}(0) = -\infty$.

Numerical analysis results show that a problem occurs if $f_0$ is greater than $10^{-3}$ and the order difference $o$ is greater than $\omega_0=22$. If these two conditions hold, we have to switch the evaluation of the integrand from \texttt{double} to \texttt{vpa}. If $f_0>10^{-3}$ and $\sigma$ is really small, then $\Re(C_2)$ is close to zero, the order difference $o=+\infty$ and we have to always switch. For values $f_0\leq 10^{-3}$, the values $f(k)$ remain close to zero for all $\Re(k)>0$ and they are correctly evaluated in \texttt{double}, so it is not necessary to switch to \texttt{vpa} and the integral will also be close to zero. Empirical value $\omega_0$ is determined by the values of the integrand $f(k)$ to be evaluated within precision \texttt{1e-5}. Such a precision together with the exponential decay of the integrand allows us to evaluate the corresponding integral \eqref{e:price} numerically within precision \texttt{1e-8} (see numerical results in Section \ref{sec:results}).

\begin{example}[Test Case 1 revisited]\label{ex:test_cases_cts2}
In Example \ref{ex:test_cases} Figure \ref{fig:test1} we presented inaccurately evaluated integrand for certain parameter values. In Table \ref{tab:tablenew} we can see the relation of the order difference $o$ and the integrand evaluation error for different values of parameter $\sigma$, all other parameters remain the same. Values of $f_0$ are evaluated both in \texttt{double} or \texttt{vpa} and their difference is in the 7th column. In Figure \ref{fig:fignew} we can see an $\text{error}_f :=\Re(f^{\texttt{double}}(k) - f^{\texttt{vpa}}(k))$ for $\Re(k)\in[0,1]$. Maxima of these differences are listed in the last column of Table \ref{tab:tablenew}. We can see that the error for $o>22$ is bigger than the desired precision \texttt{1e-5}. 

\begin{table}[ht!]
\centering
\caption{Integrand evaluation errors in Test Case 1}\label{tab:tablenew}
\begingroup
{\footnotesize
\begin{tabular}{llllllll}     
\toprule
$\sigma$    &    $o_1$   &   $o_2$   &  $o$  &  $f_0^{\texttt{double}}$  &  $f_0^{\texttt{vpa}}$ &  $|f_0^{\texttt{double}} - f_0^{\texttt{vpa}}|$ & $\max |\text{error}_f|$ \\
\midrule
0.1     & 5.061  &  -8.510 &   13.571  &    2.1369600154  &    2.1369600152  & 0.0000000001 & 0.0000000009 \\
0.05    & 5.663  &  -9.112 &   14.775  &    2.1369590951  &    2.1369590940  & 0.0000000010 & 0.0000000034 \\
0.01    & 7.061  & -10.510 &   17.571  &    2.1369583500  &    2.1369583571  & 0.0000000070 & 0.0000000887 \\
0.005   & 7.663  & -11.112 &   18.775  &    2.1369583299  &    2.1369582649  & 0.0000000650 & 0.0000003834 \\   
0.001   & 9.061  & -12.510 &   21.571  &    2.1369592635  &    2.1369581912  & 0.0000010722 & 0.0000093898 \\
0.0005  & 9.663  & -13.111 &   \color{myhalforange}{22.775}  &    2.1369775046  &    2.1369581820  & 0.0000193225 & \color{myhalfred}{0.0000378410} \\
0.0001  & 11.061 & -14.507 &   \color{myhalforange}{25.568}  &    2.1370505356  &    2.1369581747  & 0.0000923609 & \color{myhalfred}{0.0007385972} \\
0.00005 & 11.663 & -15.051 &   \color{myhalforange}{26.714}  &    2.1388775265  &    2.1369581737  & 0.0019193527 & \color{myhalfred}{0.0020950859} \\   
0.00001 & 13.061 &    -Inf &   \color{myhalforange}{Inf   }  &    2.1243052036  &    2.1369581730  & 0.0126529693 & \color{myhalfred}{0.0126529693} \\ 
\bottomrule
\end{tabular}
}
\endgroup

\end{table}

\begin{figure}
\includegraphics[width=\textwidth]{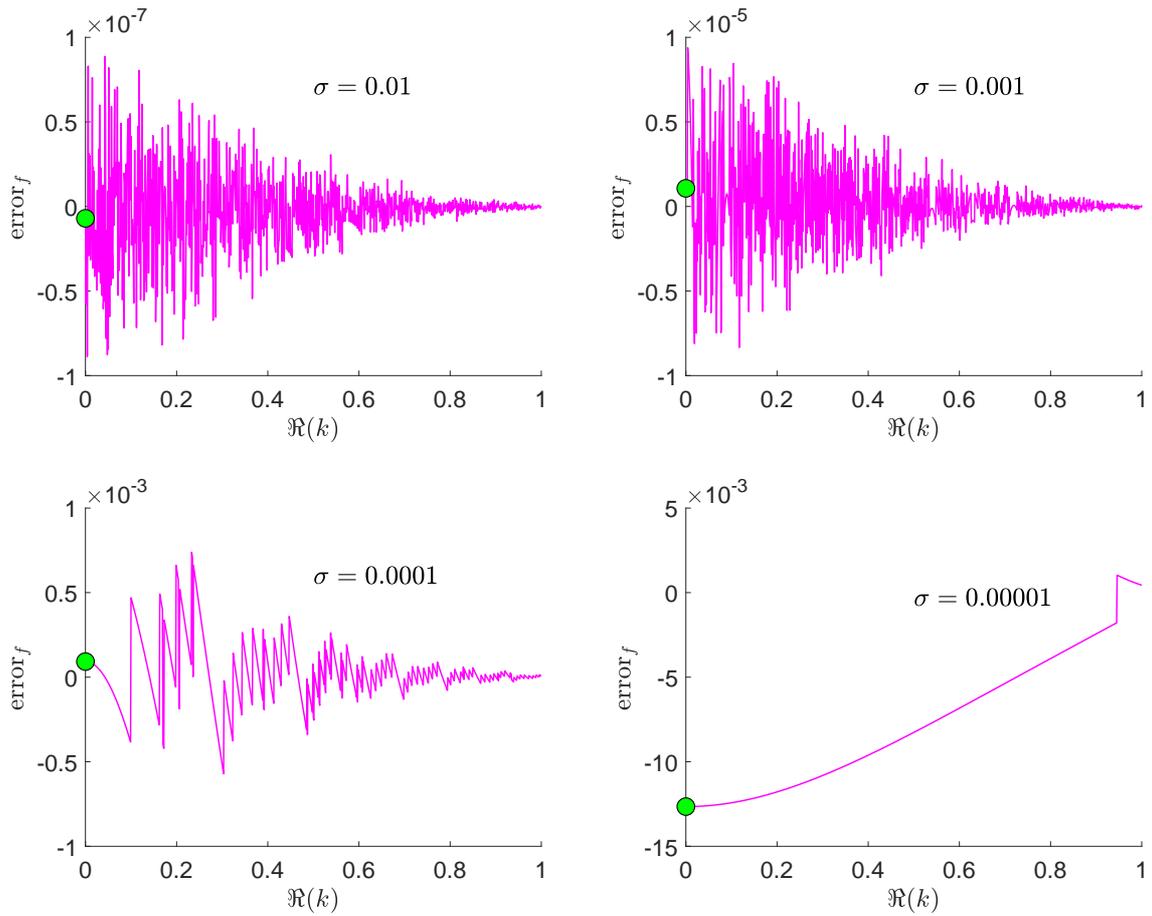} 
\caption{Integrand evaluation errors in Test Case 1 -- see Example \ref{ex:test_cases_cts2}. Values are sampled at equidistant points with step-size \texttt{1e-3}. }\label{fig:fignew}
\end{figure}

\end{example}

Further numerical analysis shows that it is not necessary to switch to \texttt{vpa} in all such cases and we can make some \emph{corrections} to the empirical value $\omega_0$. The first correction is based on the order of the term in front of the integrand. The higher the order of this term, the higher order $o$ can lead to an integrand that is sufficient to be evaluated in standard \texttt{double} arithmetic (for simplicity we write that the integrand is ``\texttt{double} sufficient''), i.e. it is not necessarily to switch to \texttt{vpa}. Since $K e^{-r\tau} / \pi < K/3$, we define
\begin{equation}\label{e:omega1}
\omega_1 := \min(5,\max(4- \log_{10}(K/3) ,0)).
\end{equation}

The second correction is based on the order of the value $f_0$. Here the influence is reverse, the higher the value of the order of $f_0$, the worse is the behaviour of the integrand, i.e. the order $o$ must be lower in order for the integrand to  be \texttt{double} sufficient. Let
\begin{equation}\label{e:omega2}
\omega_2 := \min( \log_{10}|\Re(f_0)|, 0).
\end{equation}

The overall switching regime algorithm is summarized as Algorithm \ref{alg:switch} below. It is worth to mention that the switching criteria is relatively cheap operation that requires only one additional function evaluation, namely calculation of $f_0$ in \texttt{double} arithmetic. Additional speed-up of the algorithm could be achieved in faster calculation of the order of relevant terms, i.e. instead of calculating the decimal order using the base 10 logarithm, one could for example take the exponent from the floating point representation of the number. Such a modification goes beyond the scope of this manuscript. 

\begin{algorithm}
\begin{algorithmic}
\STATE $o_1 := \log_{10}|\Re (C_1)|$;\quad $o_2 := \log_{10}|\Re (C_2)|$;\quad $o   := o_1 - o_2$;\quad  $\omega_0 = 22$;
\STATE $f_0 := f(0+i/2)$; 
\STATE \texttt{par} = \FALSE;
\IF{$f_0 > 10^{-3}$ \AND $o > \omega_0$}
\STATE $\omega_1 := \min(5,\max(4- \log_{10}(K/3) ,0))$; 
\STATE $\omega_2 := \min( \log_{10}|\Re(f_0)|, 0)$;
\IF{$o > \omega_0 + \omega_1 - \omega_2$} 
\STATE \texttt{par} = \TRUE;
\ENDIF
\ENDIF
\IF{\texttt{par}}
\STATE switch the evaluation of the integrand $f(k)$ from \texttt{double} to \texttt{vpa}
\ENDIF
\end{algorithmic}
\caption{Fast regime switching algorithm for the evaluation of the integrand $f(k)$ that is by default in \texttt{double}.}\label{alg:switch}
\end{algorithm}

\begin{figure}
\includegraphics[width=0.49\textwidth]{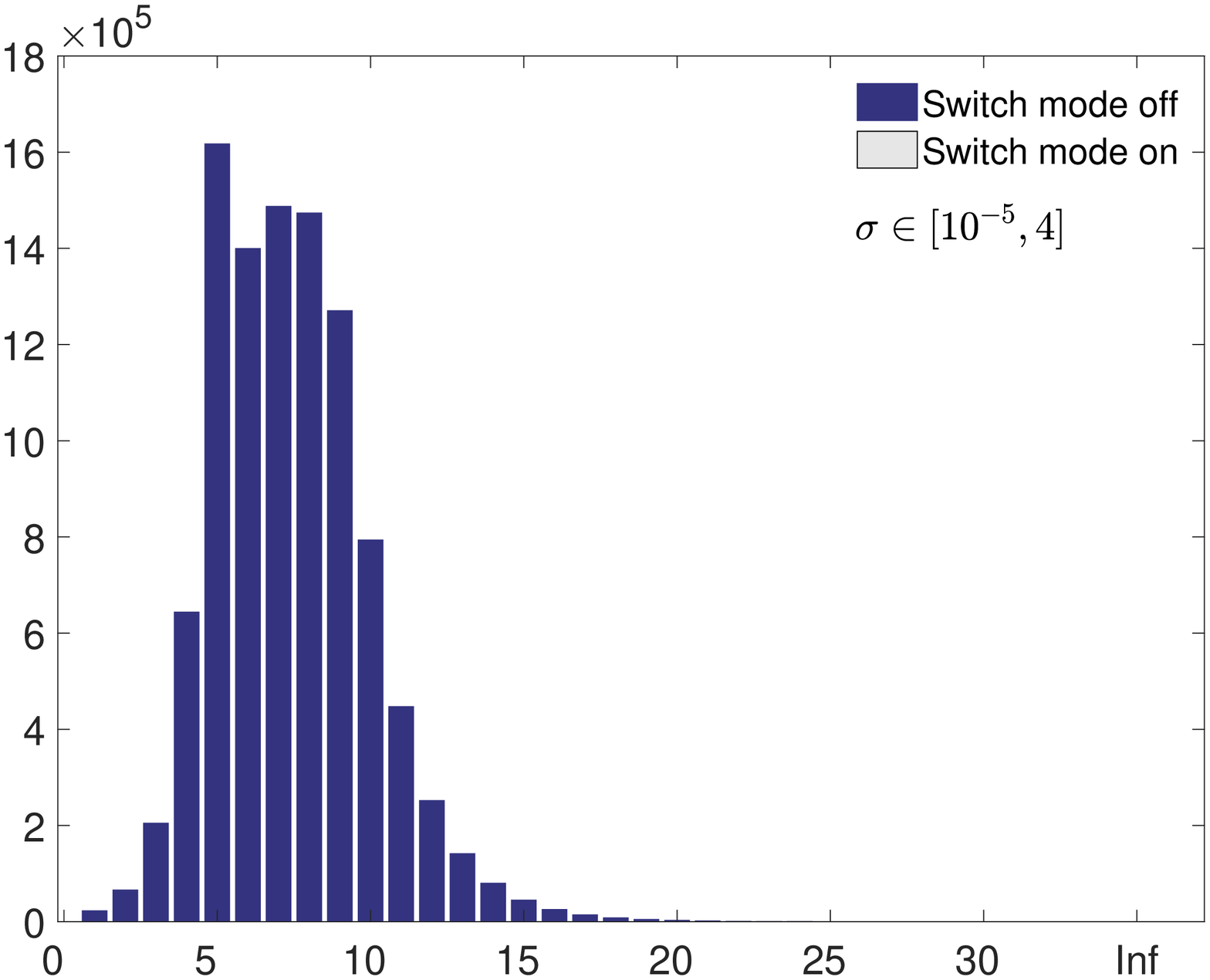}~%
\includegraphics[width=0.49\textwidth]{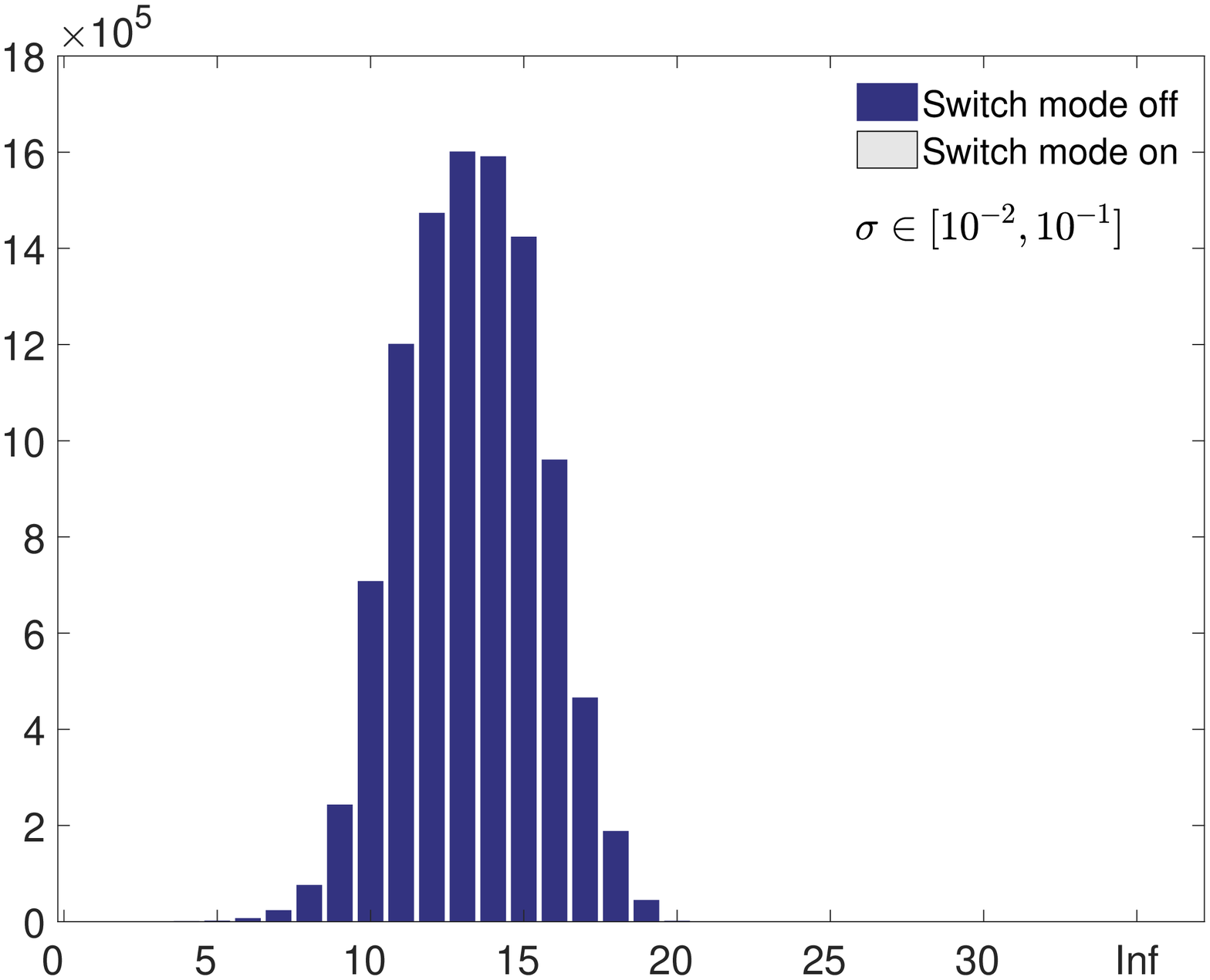}\\
\includegraphics[width=0.49\textwidth]{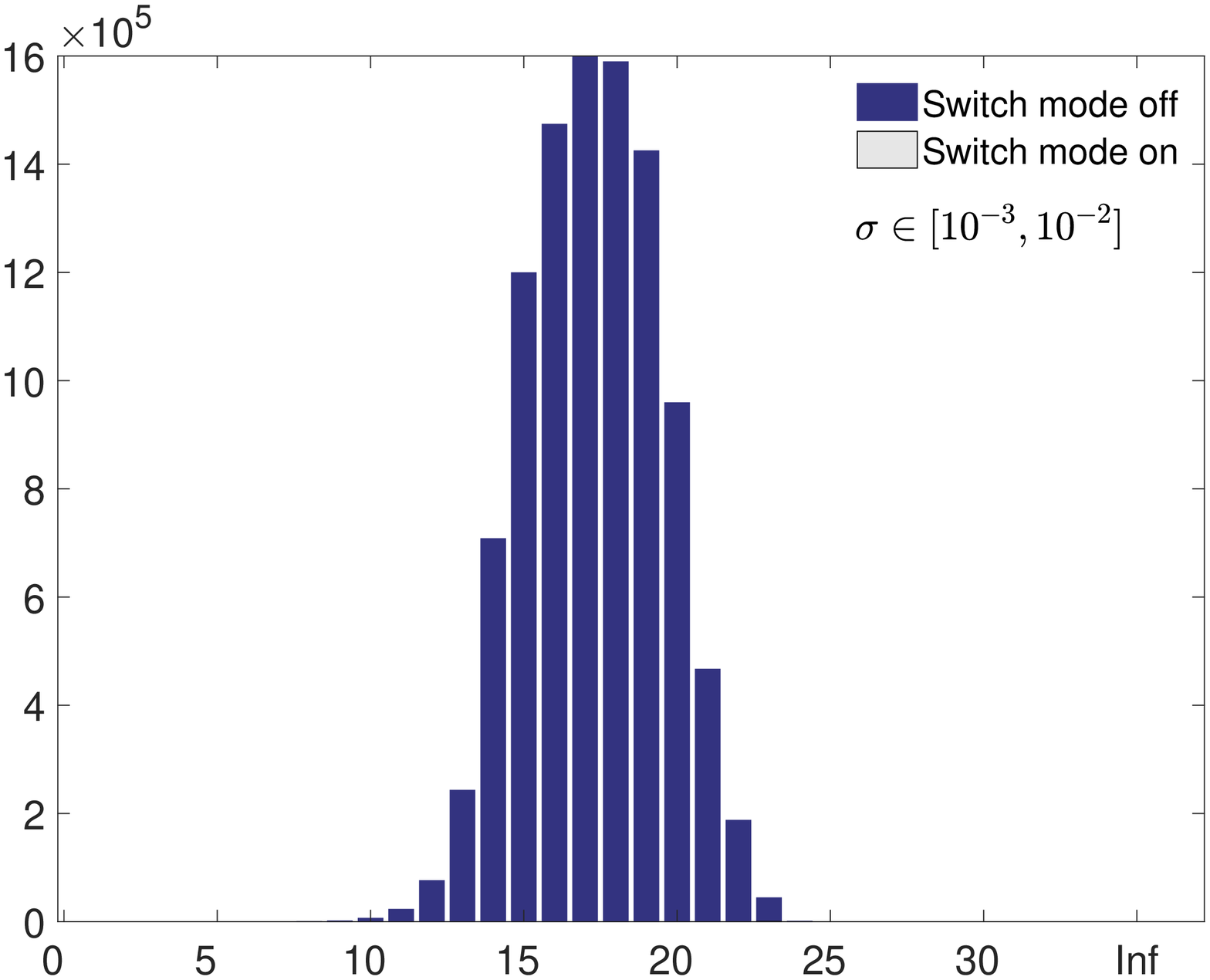}~%
\includegraphics[width=0.49\textwidth]{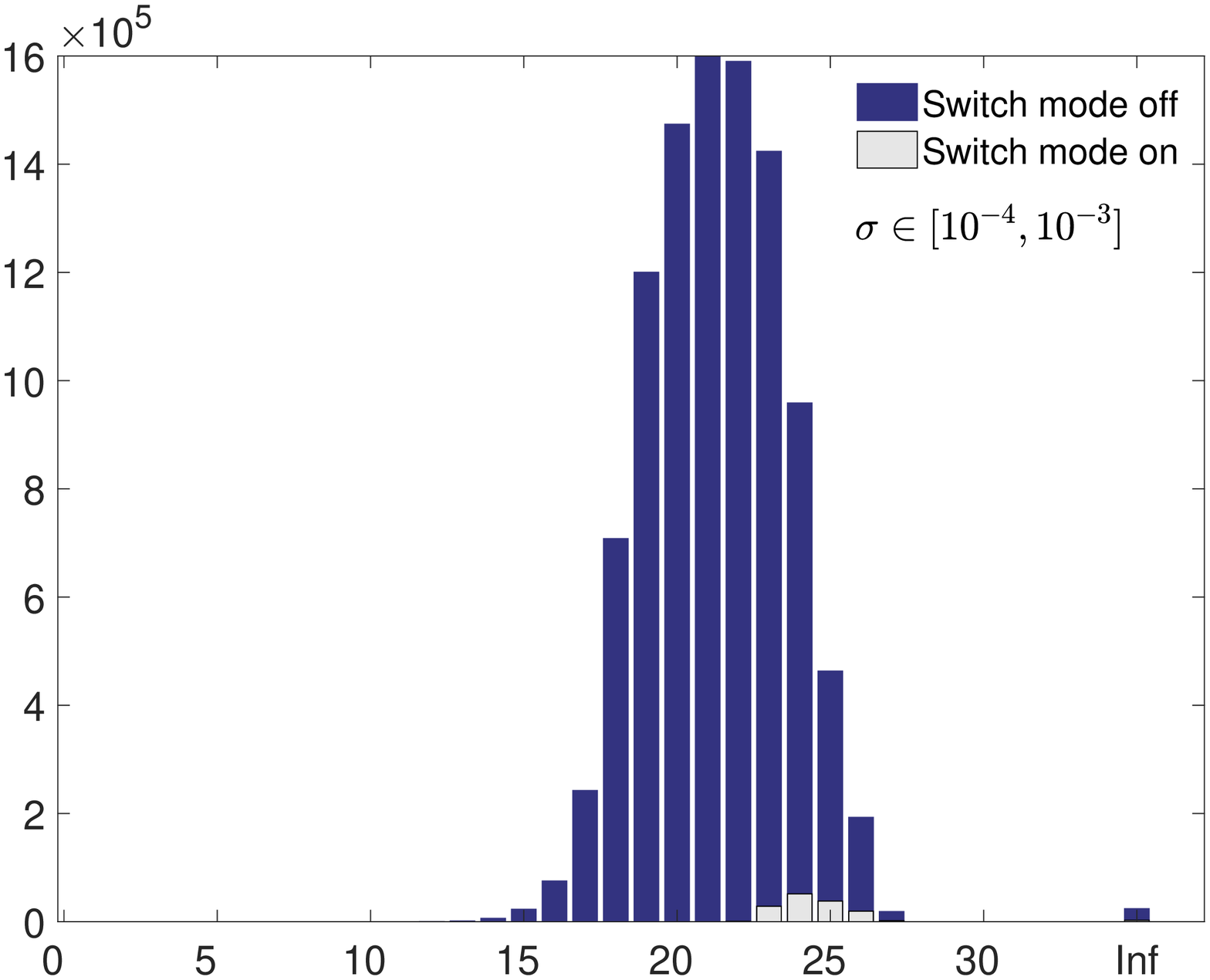}\\
\includegraphics[width=0.49\textwidth]{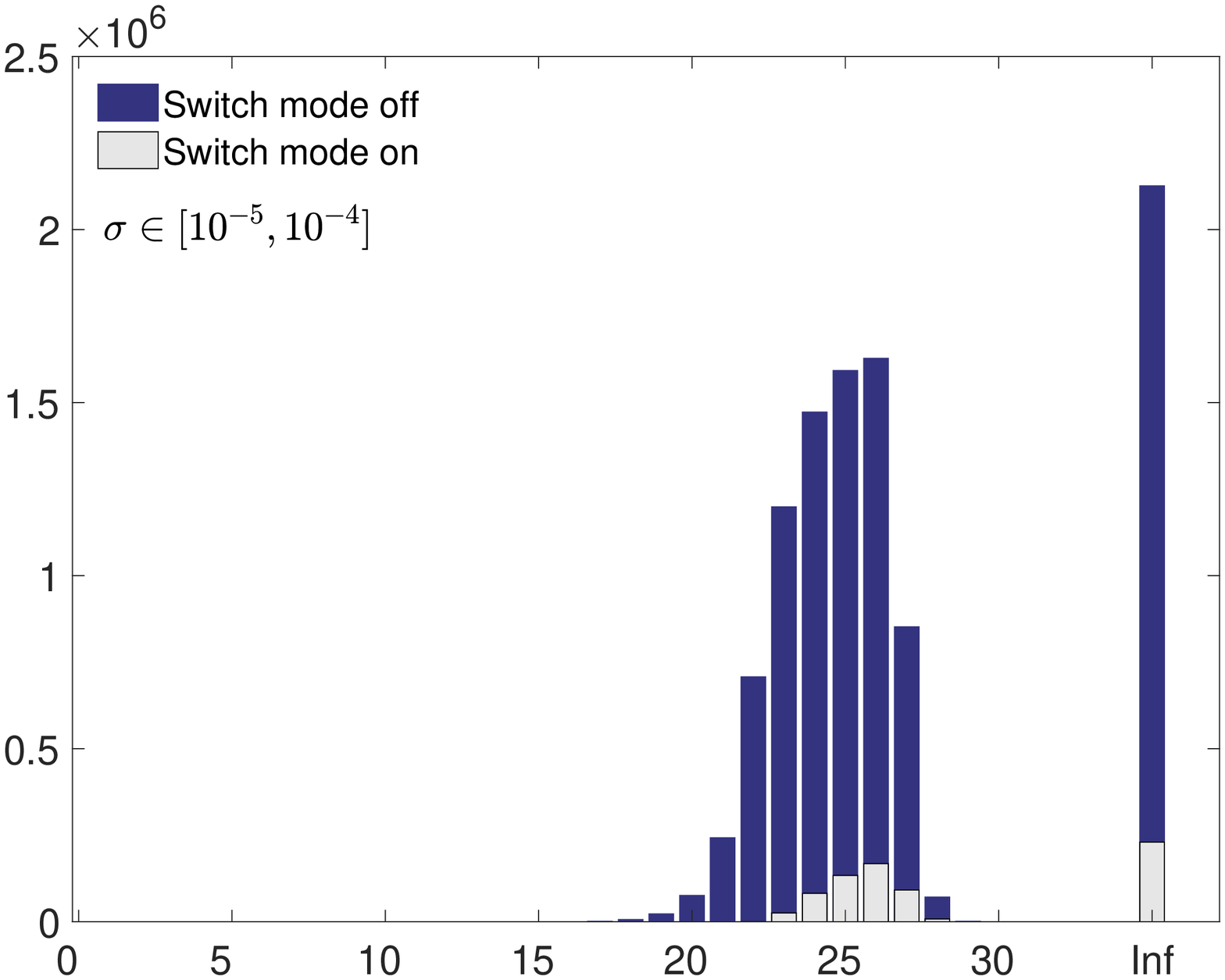}~%
\includegraphics[width=0.49\textwidth]{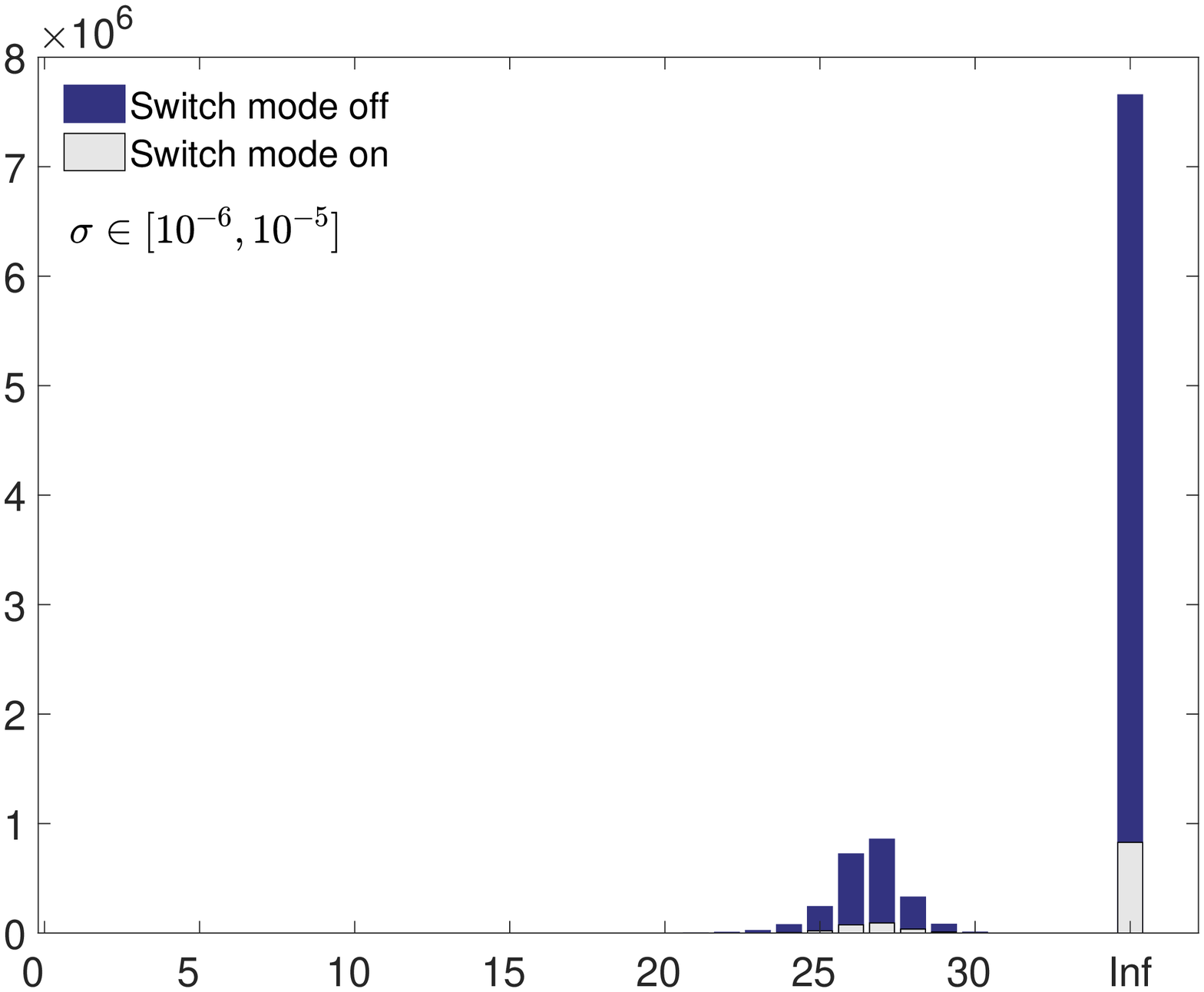}%
\caption{Influence of the parameter $\sigma$ to the switch mode. Horizontal axis shows the value of the order difference $o$, the vertical axis is the number of cases in 10 millions integral evaluations. All parameter values were generated uniformly random in considered bounds.}\label{fig:histograms}
\end{figure}

In Figure \ref{fig:histograms}, we can see the histograms for the order difference $o$. In each case, 10 millions integral evaluations were performed with parameter values taken uniformly random in considered bounds. "Switch mode on" means that in the Algorithm \ref{alg:switch} the value \texttt{par} = \textbf{true}. The lower the value of the parameter $\sigma$, the higher is the order difference $o$. However, as we can see, it is not only the value $\sigma$ that causes the problems. For many cases it is sufficient to evaluate the integrand $f(k)$ in \texttt{double} arithmetic only even if the value $\sigma$ is really low. In the last graph the value \texttt{Inf} indicate that the value $o=+\infty$, i.e. that $\Re(C_2)$ is close to zero. 

Note that histograms in Figure \ref{fig:histograms} correspond also to the values obtained in Table \ref{tab:2-1-cetnosti} and the thorough analysis of problematic cases motivated the choice of $\omega_0=22$ that can be observed also in histograms in the third row.

\begin{figure}[ht!]
\begin{center}
\includegraphics[width=0.49\textwidth]{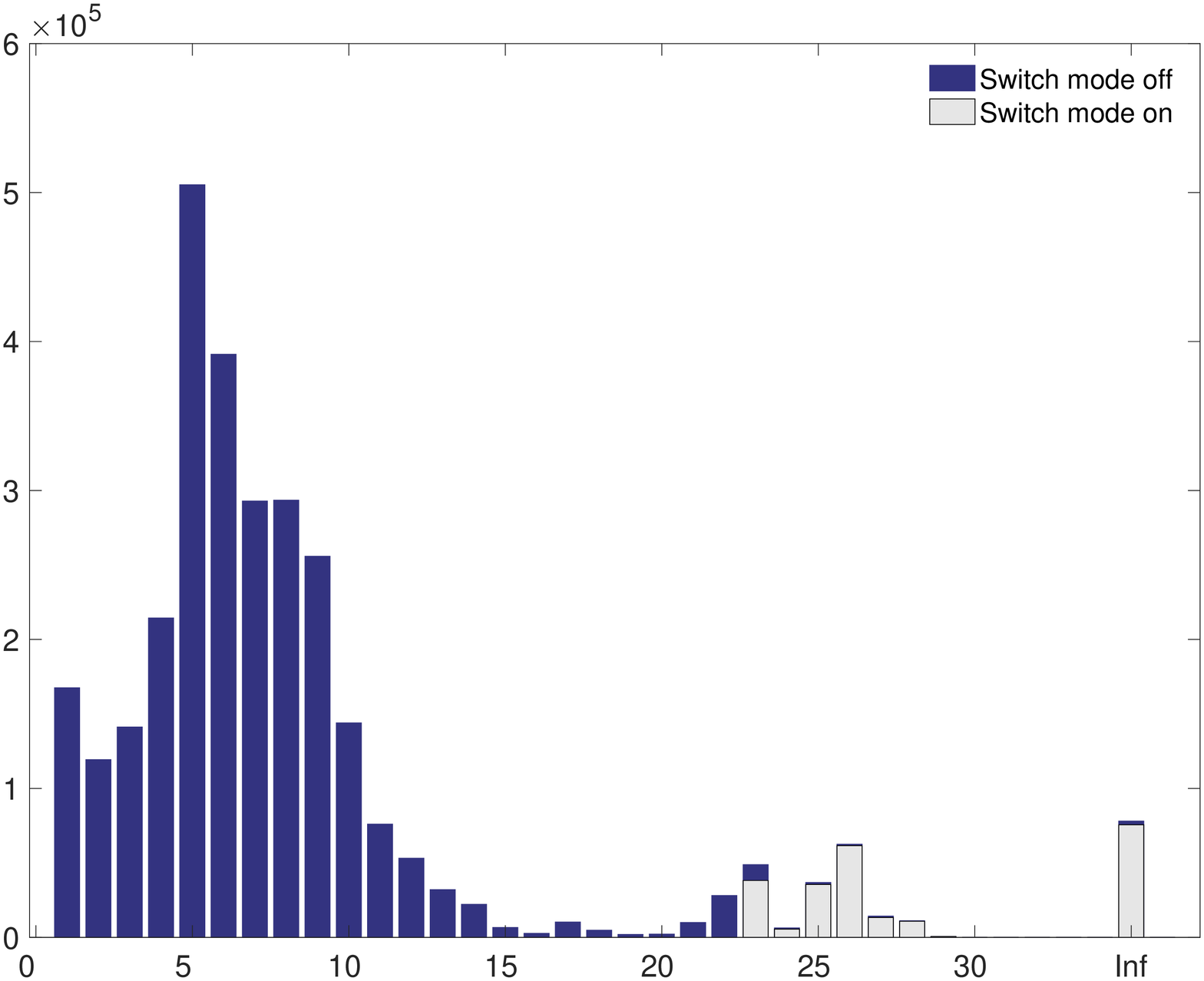}~%
\includegraphics[width=0.49\textwidth]{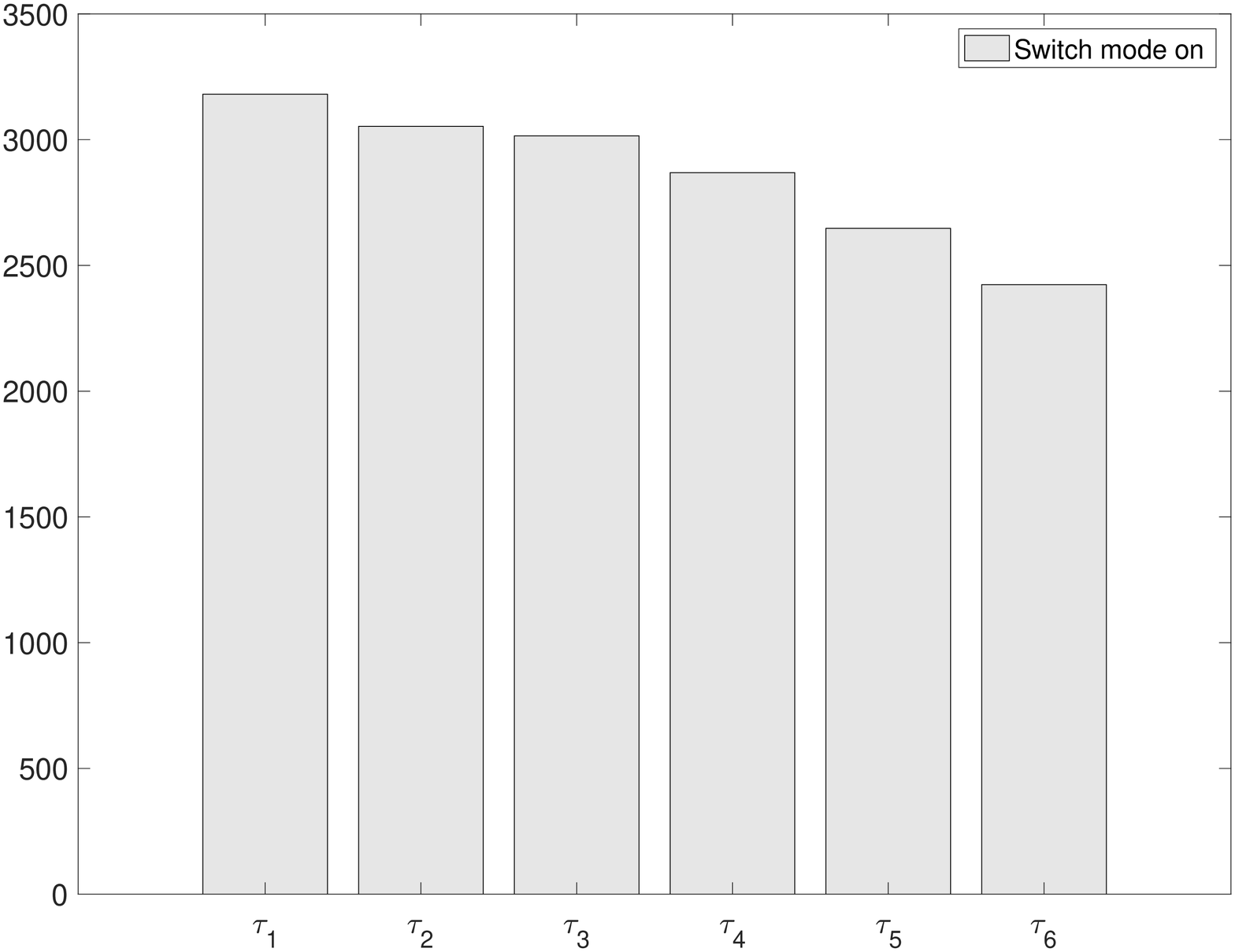}
\caption{On the left, there is a histogram showing the \texttt{double} insufficiency during 10 independent global optimization runs. Horizontal axis shows the value of the order difference $o$, the vertical axis is the number of cases in 3\,034\,000 that were performed during the optimization. On the right, an average number of \texttt{double} insufficient integrand evaluations is grouped by the time to maturity $\tau$.}\label{fig:cal}
\end{center}
\end{figure}

\begin{example}[Problematic integrand evaluations during calibration]\label{ex:cal}
Calibration of the models to real market data will be further studied below in Section \ref{ssec:cal}. In this example we give more details to the number of problematic cases that occur during calibration processes described later in Example \ref{ex:calib}. Since the random generation of parameters occur only during the global optimization phase of the calibration process, we analyse the problematic integrand evaluations during this phase only.

In particular, from ten independent global optimization runs for the AFSVJD model we gathered 3\,034\,000 different vectors of parameter values. For these vectors we tested if an integrand evaluation is \texttt{double} sufficient or not. We found out that there were 241\,307 problematic cases, i.e. almost eight percent (7.95\%). In Figure \ref{fig:cal}, we can see the \texttt{double} insufficiency distribution together with the dependence on the time to maturity (there are 6 different maturities in the data set ranging from $\tau_1 = 0.120548$ to $\tau_6 = 0.977528$). From this observation we can conclude that on average, more problematic cases occur for shorter maturities.
\end{example}

If the integrand is \texttt{double} insufficient and higher precision arithmetic is not used, the integral evaluation and hence the overall calibration process can slow down dramatically. Both accuracy and speed of the numerical quadratures will be thoroughly tested in the next section.

\clearpage
\section{Numerical quadratures and their failures}\label{sec:quadratures}

The problem in numerical integration is to approximate definite integrals of the form $\int_a^b f(x)\,dx$ by the $n$-point numerical quadrature 
\begin{equation}
\int_a^b f(x)\d x \approx Q_n[a,b] := \sum_{i=1}^n w_i f(x_i),
\end{equation}
where $w_i$ are the weights and $x_i$ are the points at which the function $f(x)$ is evaluated. An $n$-point Gaussian quadrature is a rule constructed to give an exact result for polynomials of degree $2n-1$ or less by a suitable choice of the points $x_i$ and weights $w_i$ for $i=1,\dots,n$. It can be shown \cite{Gautschi04,Press07} that the evaluation points $x_i$ are the roots of a polynomial belonging to a class of orthogonal polynomials. Gauss-Legendre $n$-point quadrature (denoted by \texttt{legendre($n$)} below) is probably the best known Gaussian quadrature with associated orthogonal polynomials being the Legendre polynomials. 

The adaptive control strategy divides the integration domain into subintervals, evaluates the integral at each region and uses an error estimate of the integral to check if a specified error tolerance is met. At regions where the function is well approximated by a polynomial, only a few function evaluations are needed, in other areas the adaptive strategy evaluates the subintervals in a recursive manner. Adaptive quadratures together with detailed error estimation techniques were reviewed by \citet{Gander00}, where adaptive Simpson quadrature (\texttt{quad}) was also introduced. This review was recently extended by \citet{Gonnet09,Gonnet12}.

Adaptive strategy is applied also in the adaptive Gauss-Lobatto quadrature with a modification by the so called Kronrod extension to add an effective error control procedure. We consider two implementations of Gauss-Lobatto, namely \texttt{adaptlob} is the original implementation by \citet{Gander00} and \texttt{quadl} is the MATLAB implementation that further improves the adaptivity. The Lobatto formula has preassigned abscissas at the end points of the interval and other nodes and weights are determined in order to obtain the highest exactness possible. The Kronrod extension is used to provide an estimate of the approximation error. If the error exceeds a specified tolerance, regions where the function is not well behaved will be divided recurrently. 

Gauss-Kronrod quadrature is an adaptive extension of the Gauss-Legendre algorithm in which the evaluation points are chosen so that an accurate approximation can be computed by re-using the information produced by the computation of a less accurate approximation. For the same set of function evaluation points, it has two quadrature rules, one higher order and an embedded one with lower order. The difference between these two approximations is used to estimate the calculation error of the integration. If the quadrature is applied to the interval $[a_k,b_k]$, then the error estimate takes the form
\begin{equation}\label{e:gkerr}
E_k = | G_n[a_k,b_k] - K_{2n+1}[a_k,b_k]|,
\end{equation}
where $G_n[a,b]$ is the $n$-point Gauss quadrature rule of degree $2n-1$, and $K_{2n+1}[a,b]$ is the $2n+1$ point Gauss-Kronrod extension of degree $3n+1$ that is used as the approximation of the integral. This is also the error estimate currently used in the MATLAB function \texttt{integral} that implements the adaptive Gauss-Kronrod(7,15) quadrature. This quadrature was chosen as the reference, since it is used de facto as an industrial standard in the latest releases of MATLAB. Initially, the interval $[a,b]$ is split into 10 equally sized subintervals. The default recurrence rule is set so that if the error estimate \eqref{e:gkerr} for an interval $[a_k,b_k]$ is greater than the prescribed tolerance, the interval $[a_k,b_k]$ is divided to halves where the process is recurrently repeated. A~reference to Gauss-Kronrod quadratures is \cite{Gautschi88}, computation of Gauss-type quadrature formulas is further studied for example in \cite{Laurie01}. It is worth to mention that other error estimates for the Gauss-Kronrod quadrature exist, some of them are of rather curious empirical character (e.g. the one used in the widely-used QUADPACK library, cf. \cite[Section 2.4]{Gonnet12}), but in principle none of the estimates can properly handle the problem of loosing significant digits precision in the inaccurate integrand evaluation. This problem leads us to another simple example showing the failure of numerical quadratures, in particular of the adaptive Gauss-Kronrod quadrature.

\begin{example}[Failure of the Gauss-Kronrod quadrature]\label{ex:gk37_example}
Let $f_0(x)$ be a piecewise constant function taking values either $1-\varepsilon$ or $1+\varepsilon$ with discontinuities exactly at the quadrature abscissas so that $\int\limits_a^b f_0(x)\d x = b-a$, see Figure \ref{fig:gk37_example} at the top where discontinuities occur at the abscissas of the Gauss-Kronrod(3,7) quadrature. The value of $\varepsilon$ is sufficiently small in order to serve as a measure of discontinuities caused by loosing significant digits. Let us consider the general Gauss-Kronrod($n$,$2n+1$) quadrature. The even nodes that partition the interval $[a,b]$ are the Gaussian points of the $G_n[a,b]$ quadrature and $f_0$ is equal to $1-\varepsilon$ at these points and $1+\varepsilon$ at the odd nodes that are the Gaussian points of the Kronrod extension $K_{2n+1}[a,b]$. We have that
\begin{align*}
G_n[a,b] &= \sum_{i=1}^{n} w_i^G f_0(x_i^G) = (1-\varepsilon) \sum_{i=1}^n w_i^G = (1-\varepsilon)(b-a) \\ 
K_{2n+1}[a,b] 
&= \sum_{i=1}^{2n+1} w_i^K f_0(x_i^K) = \sum_{i=1}^{n} w_{2i+1}^K f_0(x_{2i+1}^K) + \sum_{i=1}^{n} w_{2i}^K f_0(x_{2i}^K) \\ 
&= (1+\varepsilon) \sum_{i=1}^{n} w_{2i+1}^K + (1-\varepsilon) \sum_{i=1}^{n} w_{2i}^K 
= \sum_{i=1}^{2n+1} w_i^K + \varepsilon \left( \sum_{i=1}^{n} w_{2i+1}^K  - \sum_{i=1}^{n} w_{2i}^K \right) \\
&= b-a + \varepsilon (b-a) C_n = (b-a)(1+\varepsilon C_n)
\intertext{and the error estimate}
E_0 
&= |G_n[a,b] - K_{2n+1}[a,b]| 
= |(1-\varepsilon)(b-a) - (b-a)(1 + \varepsilon C_n)|  \\
& = \varepsilon (b-a) \tilde{C}_n.
\end{align*}
For example for $\varepsilon = 10^{-4}$, $n=7$ (i.e. for Gauss-Kronrod(7,15) quadrature) and $[a,b]=[-1,1]$ we get $E_0 \doteq 2.004652 \times 10^{-4}$ which is much bigger than the default absolute or relative tolerance \texttt{AbsTol=1e-10} and \texttt{RelTol=1e-6} respectively and hence the adaptive refinement is required. 

To mimic the inaccurate evaluation of a function caused by loosing significant digits, we can consider a piecewise constant function $f_l(x)$ taking again values either $1-\varepsilon$ or $1+\varepsilon$ with discontinuities at the quadrature abscissas in all of the $2^l$ subintervals $\left[a+\frac{i-1}{2^l}(b-a), a+\frac{i}{2^l}(b-a) \right]$, $i=1,\dots,2^l$. In Figure \ref{fig:gk37_example} at the bottom we can see a graph of function $f_2$, i.e. the case where the original interval $[a,b]$ was divided into 4 subintervals and at each of the subinterval the discontinuities occur at the seven abscissas of the Gauss-Kronrod(3,7) quadrature. To approximate the exact value $\int\limits_a^b f_l(x)\d x = b-a$ by the adaptive Gauss-Kronrod($n$,$2n+1$) quadrature, adaptive refinement must be performed at least to the $2^{l+1}$ subintervals, where the local error at the level $2^l$ is $E_0$. 

\begin{figure}[ht!]
\centering
\includegraphics[width=\textwidth]{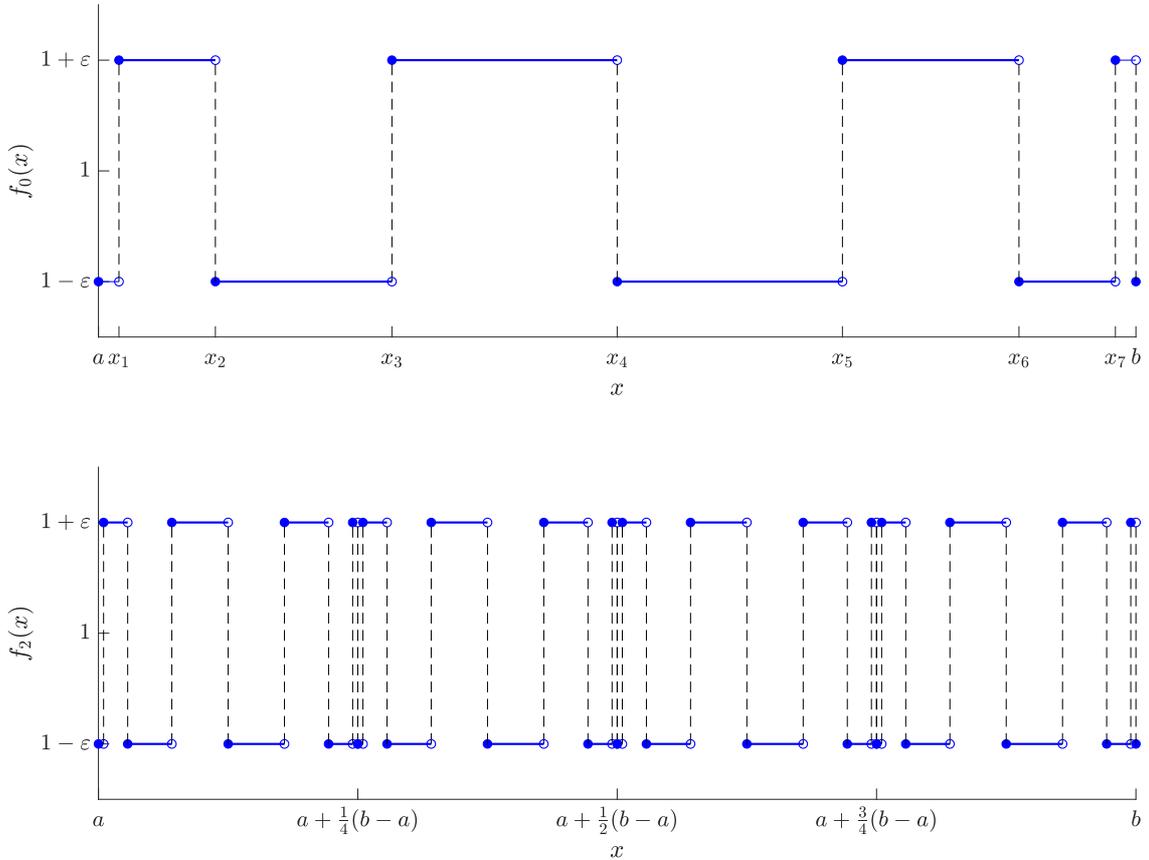}
\caption{Functions $f_0$ and $f_2$ from the Example \ref{ex:gk37_example}. Discontinuities occur at the abscissas of the Gauss-Kronrod(3,7) quadrature.}\label{fig:gk37_example}
\end{figure}

\end{example}

\begin{example}[Test Case 1 revisited]
In Example \ref{ex:test_cases} Figure \ref{fig:test1} we presented inaccurately evaluated integrand for certain parameter values that could be obtained during the calibration process. Examining the integrand for these parameter values, we can find out that in \texttt{double} evaluation there were $63\,302$ ``significant discontinuities'' of $f(k)$ for $\Re(k)\in[0,20]$ with jump sizes up to $8.6846\times 10^{-4}$ (average jump size $1.9082\times 10^{-6}$). To satisfy the local Gauss-Kronrod error estimate to be within the relative error $10^{-6}$, adaptive refinement had to be done down to the interval length $ 5\times 10^{-8}$. In the language of the previous Example \ref{ex:gk37_example}, the value of $l$ can grow up to 12--16. Not only is such a numerical integration time demanding (one integral can be calculated several seconds, for exact values see results in Section \ref{sec:results}), the obtained integral value is far from the precise value, the error is $2.4679\times 10^{-5}$ although the adaptive quadrature ends with errors within the default tolerances \texttt{AbsTol=1e-10} and \texttt{RelTol=1e-6}.
\end{example}

From the review of the adaptive quadratures written by \citet{Gonnet09,Gonnet12} it is obvious that although many other error estimates exist, none of them can properly handle the inaccuracy of the function evaluation caused by loosing significant digits. Almost all existing error estimates assume that the integrand is sufficiently smooth. The only exception is probably the \texttt{CADRE} estimate by de Boor (see \cite[Section 3.1]{Gonnet12}) that allows to detect the low frequency discontinuities. Although the integrand in our case is theoretically sufficiently smooth, for some parameter values its evaluation is numerically \texttt{double} insufficient. The only way how to calculate the integral of such an integrand is to use the high precision arithmetic. Much simpler integrand is for example studied by \citet{Gautschi16acm} who introduces a non-standard Gauss-Hermite quadrature with empirical error estimates and implementation \cite{Gautschi17opq} using the MATLAB \texttt{vpa}.

In the numerical comparison provided in the next section, \texttt{integral (vpa)} denotes the case when the function \texttt{integral} is applied to \texttt{vpa} evaluated integrand, i.e. the quadrature uses values that are evaluated in \texttt{vpa} and converted to \texttt{double}. In the MATLAB release R2016b, a new function \texttt{vpaintegral} was introduced. This function implements a semi-symbolic quadrature that is supposed to be the high-precision numerical integration. Although the way how this quadrature obtains the numerical result differs from all other studied quadratures, we also consider it in our comparisons. For convenience we also add the very simple trapezoidal rule (\texttt{trapz}) that together with \texttt{legendre} is not adaptive. To understand how the adaptivity works in all considered quadratures is crucial to realize what happens in a numerical quadrature if the integrand is inaccurately evaluated as we saw for example in Figures \ref{fig:test1} and \ref{fig:test2}.

\section{Numerical results}\label{sec:results}

In this section we compare numerical quadratures behaviour in problematic cases and with the implementation of the switching regime Algorithm \ref{alg:switch} described in Section \ref{sec:integrand}.

We measure an average calculation time (reference PC with 1x quad-core Intel i7-4770K 3.5 GHz CPU and 16 GB RAM) and number of integrand evaluations (\texttt{fevals}). Interesting results are highlighted. In all quadratures we set the tolerances to satisfy requirements mentioned in Section \ref{sec:integrand}, i.e. we set absolute tolerance to be \texttt{AbsTol=1e-10} and relative tolerance \texttt{RelTol=1e-6}.

\subsection{Quadratures behaviour in problematic cases}

In Table \ref{tab:table1} we can see quadratures comparison for numerical integration with parameter values taken from test case 1, see Figure \ref{fig:test1}. Although the global view to the integrand does not indicate any potential problem, the contrary is true. We compare several quadratures for \texttt{double} evaluated integrand with reference value \texttt{integral (vpa)}, i.e. \texttt{integral} function applied to the \texttt{vpa} evaluated integrand. In fact, the whole integrand $f(k)$ is evaluated using \texttt{vpa} and before plugging it into the \texttt{integral} function, it is converted into \texttt{double}. For convenience, we also provide the simple \texttt{trapz} rule with step sizes 0.01 and 0.001 and \texttt{quadl} all for \texttt{double} or \texttt{vpa} evaluated integrand.

In the case for $\sigma=0.001$, apart from \texttt{quadl} errors for all other quadratures are negligible, as well as computation times (without \texttt{vpa}) and \texttt{fevals} are small. Also \texttt{quadl (vpa)} give us accurate result (in given tolerance \texttt{1e-8}). From the \texttt{integral} and \texttt{integral (vpa)} comparison we can see that even in this example the value differs. Computation time for \texttt{vpa} is ca 100 times higher which goes against our requirement to calculate the integral quickly. We can see that using \texttt{vpa} evaluated integrand in this case is not necessary for \texttt{integral}, but it is necessary for \texttt{quadl}. Computation times are especially huge for \texttt{trapz (vpa)}.

A~problematic case occurs if we change $\sigma$ to $0.0001$. Errors for all quadratures increase due to the inaccurately evaluated integrands. Computation time for \texttt{integral} is larger than in the previous case and it is actually larger than for \texttt{integral (vpa)}, which is caused by enormous \texttt{fevals} increase. Such a huge increase is caused by the adaptivity of the quadrature applied to the inaccurately evaluated integrand. 

Further decrease of $\sigma$ to $0.00001$ causes further error increase for all quadratures. Although one could expect larger computation time and \texttt{fevals}, these values for \texttt{integral} are much smaller than in the previous case, but error is much larger. Even if the computation time is low, there can be still a problematic case with high precision error and one can recognize it by large \texttt{fevals}. This is the reason why one of our criteria for problematic cases is  \texttt{fevals} larger than \texttt{1e4}.

From all three trials we learned that if we use \texttt{vpa}, \texttt{fevals} remain the same and as well as the computation times are very similar. This is an important indicator that the problematic cases are caused by the inaccurately evaluated integrands.

\begin{table}
\centering
\caption{Quadratures comparison for Test Case 1}\label{tab:table1}
\begingroup
{\footnotesize
\begin{tabular}{lllcrr}      
\toprule
{\bf method} (case $\sigma=0.001$) 
                         & {\bf value} & {\bf error}  & {\bf time [s]} & {\bf fevals} \\
\midrule
{\tt integral}           & 0.77681485  & 0.00000006  &   0.056        & 300  \\
{\tt integral (vpa)}     & 0.77681478  & 0.00000000  &   {\bf 4.948}  & 270  \\
{\tt legendre (128)}     & 0.77680533  & 0.00000944  &   0.006        &  64  \\
{\tt legendre (256)}     & 0.77682245  & 0.00000767  &   0.013        & 128  \\
{\tt quad}               & 0.77681499  & 0.00000020  &   0.092        &  1694  \\
{\tt quadl}              & 0.94027672  & {\color{myfullred}0.16346193}  &   0.150        &  2018  \\
{\tt quadl (vpa)}        & 0.77681478  & 0.00000000  &   4.110        &  146  \\
{\tt adaptlob}           & 0.77681496  & 0.00000017  &   1.485        &  30734  \\
{\tt trapz (0.01)}       & 0.77681499  & 0.00000021  &   0.014        &   1  \\
{\tt trapz (0.001)}      & 0.77681499  & 0.00000020  &   0.027        &   1  \\
{\tt trapz (vpa, 0.01)}  & 0.77681478  & 0.00000000  &  14.287        &   1  \\
{\tt trapz (vpa, 0.001)} & 0.77681478  & 0.00000000  & 301.293        &   1  \\
\midrule
{\bf method} (case $\sigma=0.0001$) 
                         & {\bf value} & {\bf error} &  {\bf time [s]} & {\bf fevals} \\
\midrule
{\tt integral}           & 0.77683946  & {\color{myhalfred}0.00002467} & {\color{myfullorange}6.761} & {\color{myfullorange}360780} \\
{\tt integral (vpa)}     & 0.77681478  & 0.00000000  &    {\bf 4.872} &     270 \\
{\tt legendre (128)}     & 0.77695229  & 0.00013750  &    0.006       &      64 \\
{\tt legendre (256)}     & 0.77682884  & 0.00001405  &    0.012       &     128 \\
{\tt quad}               & 0.77688362  & 0.00006883  &    0.218       &      5006 \\
{\tt quadl}              & 0.75020018  & {\color{myfullred}0.02661460}  &    0.146       &      2030 \\
{\tt quadl (vpa)}        & 0.77681478  & 0.00000000  &    4.101       &      146 \\
{\tt adaptlob}           & 0.77684007  & 0.00002528  &    2.278       &     47324\\
{\tt trapz (0.01)}       & 0.77684249  & 0.00002770  &    0.010       &       1 \\
{\tt trapz (0.001)}      & 0.77683949  & 0.00002470  &    0.029       &       1 \\
{\tt trapz (vpa, 0.01)}  & 0.77681478  & 0.00000000  &   13.743       &       1 \\
{\tt trapz (vpa, 0.001)} & 0.77681478  & 0.00000000  &  288.313       &       1 \\
\midrule
{\bf method} (case $\sigma=0.00001$)    
                         & {\bf value} & {\bf error} &  {\bf time [s]} & {\bf fevals} \\
\midrule
{\tt integral}           & 0.76823548  & {\color{myfullred}0.00857930} & {\color{myhalforange}0.655} & {\color{myhalforange}28260}    \\
{\tt integral (vpa)}     & 0.77681478  & 0.00000000  &   {\bf 5.004}  &   270 \\
{\tt legendre (128)}     & 0.76877555  & 0.00803923  &   0.007        &    64 \\
{\tt legendre (256)}     & 0.76810697  & 0.00870781  &   0.015        &   128 \\
{\tt quad}               & 0.76823622  & 0.00857856  &   0.099        &   2022 \\
{\tt quadl}              & 0.76846917  & 0.00834561  &   0.149        &   2018 \\
{\tt quadl (vpa)}        & 0.77681478  & 0.00000000  &   4.084        &   146 \\
{\tt adaptlob}           & 0.76823478  & 0.00858000  &   0.216        &   4232 \\
{\tt trapz (0.01)}       & 0.76828142  & 0.00853336  &   0.013        &     1 \\
{\tt trapz (0.001)}      & 0.76823697  & 0.00857781  &   0.030        &     1 \\
{\tt trapz (vpa, 0.01)}  & 0.77681478  & 0.00000000  &  14.618        &     1 \\
{\tt trapz (vpa, 0.001)} & 0.77681478  & 0.00000000  & 310.098        &     1 \\
\bottomrule
\end{tabular}
}
\endgroup

\end{table}

\subsection{Results for switching regime}

In this section we apply the switching regime Algorithm \ref{alg:switch} to the function \texttt{integral}, i.e. for problematic cases the quadrature is applied to the \texttt{vpa} evaluated integrand and in other cases standard \texttt{double} arithmetic is used. The comparison is summarized in Table \ref{tab:table4}. As we can see, the error increases for smaller values of parameter $\sigma$ which is caused by the increasing order difference~$o$. 

In Test Case 1 we can observe that a problem starts to emerge even for larger value of $\sigma=0.001$, see the first double-row where the \texttt{fevals} difference is already 30, although this case does not fall into the problematic category, the error is low and it is actually not necessary to switch to \texttt{vpa}. In other cases the difference in \texttt{fevals} is huge, apart from the last case when the error was large. In Test Case 2 for $\sigma=0.00001$ the error is negligible in given tolerance \texttt{1e-8}, however this must be interpreted as a coincidence, because the \texttt{fevals} for \texttt{integral} without \texttt{vpa} is huge.

The obtained values confirm the choice of the value $\omega_0=22$ in the switching regime Algorithm \ref{alg:switch} as well as they correspond to the choice of correction terms $\omega_1$ and $\omega_2$ whose influence was described in detail in Section \ref{sec:integrand}. It is useful to remind that the Algorithm \ref{alg:switch} is designed to be fast, i.e. for some rare parameter values it can happen that the algorithm switches to the \texttt{vpa} evaluation of the integrand even if it is \texttt{double} sufficient, however the algorithm should not forget to switch any problematic case which makes the algorithm fast and reliable. 

\begin{table}
\centering
\caption{Using switching regime \texttt{integral} in Test Case 1 and 2}\label{tab:table4}
\begingroup
{\footnotesize
\begin{tabular}{lllccrrc}      
\toprule
$\sigma$ & $o$   & $\omega_1$ & {\bf par} & \multicolumn{3}{c}{{\tt integral}\ \ / \ \ {\tt integral (vpa)}}& {\bf error} \\  
\cmidrule(l){5-7}
         & $f_0$ & $\omega_2$ &           &  {\bf value} & {\bf time} & {\bf fevals}  & \\ 
\midrule
\multicolumn{5}{l}{Test Case 1 (see Figure \ref{fig:test1})} \\
\midrule
0.001    & 21.572   & --     & {\bf false} & 0.77681485 & 0.056 & 300   &  0.00000006 \\
         & 2.137    & --     &             & 0.77681478 & 4.948 & 270   &             \\
0.0005   & 22.775   & 0.681  & {\bf true}  & 0.77681452 & {\color{myhalforange}4.241} & {\color{myhalforange}115320}&  0.00000025 \\
         & 2.137    & 0      &             & 0.77681478 & 2.244 & 270   &             \\
0.0001   & 25.568   & 0.681  & {\bf true}  & 0.77683946 & {\color{myfullorange}6.761}& {\color{myfullorange}360780}&  {\color{myhalfred}0.00002467} \\
         & 2.1371   & 0      &             & 0.77681478 & 4.872 & 270   &             \\
0.00005  & 26.715   & 0.681  & {\bf true}  & 0.77684142 & {\color{myfullorange}4.198} & {\color{myfullorange}217170}&  {\color{myhalfred}0.00002663} \\
         & 2.1389   & 0      &             & 0.77681478 & 5.056 & 270   &             \\
0.00001  & {\tt Inf}& 0.681  & {\bf true}  & 0.76823548 & {\color{myhalforange}0.655} & {\color{myhalforange}28260} &  {\color{myfullred}0.00857930} \\
         & 2.1243   & 0      &             & 0.77681478 & 5.004 & 270   &             \\
0.000001 & {\tt Inf}& 0.681  & {\bf true}  & 0.77674994 & 0.068 & 750   &  {\color{myhalfred}0.00006484} \\
         & 2.1243   & 0      &             & 0.77681478 & 5.043 & 270   &             \\[3pt]
\midrule
\multicolumn{5}{l}{Test Case 2 (see Figure \ref{fig:test2})} \\
\midrule
0.0001   & 20.447   & --     & {\bf false} & 0.00695940 & 0.066 & 900   &  0.00000000 \\
         & 0.12193  & --     &             & 0.00695940 & 16.198& 900   &             \\
0.00005  & 21.651   & --     & {\bf false} & 0.00695940 & 0.066 & 900   &  0.00000000 \\
         & 0.12193  & --     &             & 0.00695940 & 16.105& 900   &             \\
0.00001  & 24.447   & 0.681  & {\bf true}  & 0.00695940 & {\color{myfullorange}5.901} & 
{\color{myfullorange}312240} &  0.00000000 \\
         & 0.12193  & -0.914 &             & 0.00695940 & 16.118& 900   &             \\
0.000001 & 28.555   & 0.681  & {\bf true}  & 0.00695451 & {\color{myfullorange}10.316} & {\color{myfullorange}561750} &  {\color{myhalfred}0.00000489} \\
         & 0.12193  & -0.914 &             & 0.00695940 & 16.132 & 900   &             \\
\bottomrule
\end{tabular}
}
\endgroup
\end{table}

\subsection{Optimal switching regime quadratures}

So far the whole calculation of $f(k)$ was either performed in \texttt{double} or \texttt{vpa} arithmetic. Since the \texttt{vpa} is rather time consuming, additional speed-up (denoted by \texttt{opt} in tables below) can be achieved by further implementation improvements. Since the most problematic part of $f(k)$ is the term $C(k,\tau)$, we can calculate precisely (using \texttt{vpa}) only this term $C$ and then convert it to \texttt{double}. The remaining terms are \texttt{double} sufficient.

Let us now compare numerical quadrature in cases when we have to switch to the \texttt{vpa} evaluation of the integrand. In Tables \ref{tab:table5} and \ref{tab:table6} we can see all such cases that we already met in Table \ref{tab:table4} now for different numerical quadratures. Using check mark {\color{mygreen}\checkmark} we indicate that the quadrature gave us the accurate result (in given tolerance \texttt{1e-8}) with respect to the reference value obtained by the \texttt{integral (vpa)}, i.e. using the \texttt{double} precision arithmetic function \texttt{integral} applied to the precisely evaluated integrand (i.e. to integrand evaluated using \texttt{vpa} and converted back to \texttt{double}). The last row in each block shows the results for the semi-symbolic quadrature \texttt{vpaintegral} that is supposed to be the implementation if the high-precision numerical integration. Since the way how \texttt{vpaintegral} evaluates the integrand differs from all other numerical quadratures, the number of \texttt{fevals} in Tables \ref{tab:table5} and \ref{tab:table6} is omitted for this quadrature.

We can see that for all numerical quadratures computational time is proportional to \texttt{fevals}. In all cases we also got a desirable accuracy with the only exception \texttt{quad} in Test Case 2. From this point of view, results for all numerical quadratures are comparable. On the other hand, the semi-symbolic quadrature \texttt{vpaintegral} seems to be unusable and although its idea seems to be promising for the future, in its current state we cannot recommend its usage. First of all, even for the same tolerance setting (\texttt{AbsTol=1e-10, RelTol=1e-6}) it does not satisfy our precision criteria. Moreover, we can observe that the error of \texttt{vpaintegral} somehow correlates to the error of \texttt{integral} without \texttt{vpa}, i.e. to only \texttt{double} evaluated integrand and \texttt{double} calculated integral. Such a correlation indicate probably a systematic misbehaviour of \texttt{vpaintegral}.

To sum up the comparison results for other numerical quadratures, in all our tests the \texttt{quadl (opt)} quadrature performed best, since it needed fewest \texttt{fevals} and hence it was fastest in all our examples. However, as we saw in Table \ref{tab:table1}, it can have serious problems in the cases that are not switched to \texttt{vpa}, i.e. in cases where integrand is \texttt{double} sufficient. Since all other quadratures are time comparable as well and the differences are not big, we can conclude that the reference Gauss-Kronrod quadrature represented by the \texttt{integral (opt)} was chosen reasonably.

\begin{table}
\centering
\caption{Using optimal switching regime quadratures in Test Case 1}\label{tab:table5}
\begingroup
{\footnotesize
\begin{tabular}{lcccrr}      
\toprule
{\bf method}            & &   {\bf value}  & {\bf error}  &  {\bf time [s]} & {\bf fevals} \\
\midrule
\multicolumn{6}{l}{($\sigma=0.0005$, $o=22.775$, $f_0=2.137$, $\omega_1=0.681$, $\omega_2=0$, {\bf par = true})}  \\
\midrule
{\tt integral}          &                          & 0.77681452  &  {\color{myhalfred}0.00000025} &  {\color{myhalforange}2.244}  & {\color{myfullorange}115380} \\
{\tt integral (vpa)}    &                          & 0.77681478  &  0.00000000 &   4.920 &     270 \\
{\tt integral (opt)}    &\color{mygreen}\checkmark & 0.77681478  &  0.00000000 &   4.157 &     270 \\
{\tt quad (opt)}        &\color{mygreen}\checkmark & 0.77681478  &  0.00000000 &   6.206 &     304 \\
{\tt quadl (opt)}       &\color{mygreen}\checkmark & 0.77681478  &  0.00000000 &   3.520 &     146 \\
{\tt adaptlob (opt)}    &\color{mygreen}\checkmark & 0.77681478  &  0.00000000 &   3.668 &     164 \\
{\tt vpaintegral}       &                          & 0.77681305  &  {\color{myhalfred}0.00000172} &   7.694 &      - \\
\midrule

\multicolumn{6}{l}{($\sigma=0.0001$, $o=25.568$, $f_0=2.1371$, $\omega_1=0.681$, $\omega_2=0$, {\bf par = true})}  \\
\midrule
{\tt integral}          &                          & 0.77683946  &  {\color{myhalfred}0.00002467} &  {\color{myfullorange}6.761}  & {\color{myfullorange}360780} \\
{\tt integral (vpa)}    &                          & 0.77681478  &  0.00000000 &   4.872 &     270 \\
{\tt integral (opt)}    &\color{mygreen}\checkmark & 0.77681478  &  0.00000000 &   4.112 &     270 \\
{\tt quad (opt)}        &\color{mygreen}\checkmark & 0.77681478  &  0.00000000 &   6.270 &     304 \\
{\tt quadl (opt)}       &\color{mygreen}\checkmark & 0.77681478  &  0.00000000 &   3.520 &     146 \\
{\tt adaptlob (opt)}    &\color{mygreen}\checkmark & 0.77681478  &  0.00000000 &   3.467 &     164 \\
{\tt vpaintegral}       &                          & 0.77685918  &  {\color{myhalfred}0.00004440} &   31.627 &     - \\
\midrule
\multicolumn{6}{l}{($\sigma=0.00005$, $o=26.715$, $f_0=2.1389$,  $\omega_1=0.681$, $\omega_2=0$, {\bf par = true})}  \\
\midrule
{\tt integral}          &                          & 0.77684142  &  {\color{myhalfred}0.00002663} &  {\color{myfullorange}4.198}   & {\color{myfullorange}217170} \\
{\tt integral (vpa)}    &                          & 0.77681478  &  0.00000000 &  5.056 &     270 \\
{\tt integral (opt)}    &\color{mygreen}\checkmark & 0.77681478  &  0.00000000 &  4.165 &     270 \\
{\tt quad (opt)}        &\color{mygreen}\checkmark & 0.77681478  &  0.00000000 &  6.317 &     304 \\
{\tt quadl (opt)}       &\color{mygreen}\checkmark & 0.77681478  &  0.00000000 &  3.615 &     146 \\
{\tt adaptlob (opt)}    &\color{mygreen}\checkmark & 0.77681478  &  0.00000000 &  3.480 &     164 \\
{\tt vpaintegral}       &                         & 0.77682241  &  {\color{myhalfred}0.00000762} &  31.843 &      -
 \\
\midrule
\multicolumn{6}{l}{($\sigma=0.00001$, $o=\mathrm{{\tt Inf}}$, $f_0=2.1243$, $\omega_1=0.681$, $\omega_2=0$, {\bf par = true})}  \\
\midrule
{\tt integral}          &                          & 0.76823548  &  {\color{myfullred}0.00857930} &  0.655   & {\color{myhalforange}28260} \\
{\tt integral (vpa)}    &                          & 0.77681478  &  0.00000000 &  5.004 &    270 \\
{\tt integral (opt)}    &\color{mygreen}\checkmark & 0.77681478  &  0.00000000 &  4.108 &    270 \\
{\tt quad (opt)}        &\color{mygreen}\checkmark & 0.77681478  &  0.00000000 &  6.391 &    304 \\
{\tt quadl (opt)}       &\color{mygreen}\checkmark & 0.77681478  &  0.00000000 &  3.594 &    146 \\
{\tt adaptlob (opt)}    &\color{mygreen}\checkmark & 0.77681478  &  0.00000000 &  3.582 &    164 \\
{\tt  vpaintegral}      &                         & 0.76839158  &   {\color{myfullred}0.0084231} &  13.270&     -  \\
\midrule
\multicolumn{6}{l}{($\sigma=0.000001$, $o=\mathrm{{\tt Inf}}$, $f_0=2.1243$, $\omega_1=0.681$, $\omega_2=0$, {\bf par = true})}  \\
\midrule
\vrule height3pt width0pt 
{\tt integral}          &                          & 0.77674994  &  {\color{myhalfred}0.00006484} &  0.068 &  750 \\
{\tt integral (vpa)}    &                          & 0.77681478  &  0.00000000 &  5.043 &  270 \\
{\tt integral (opt)}    &\color{mygreen}\checkmark & 0.77681478  &  0.00000000 &  4.086 &  270 \\
{\tt quad (opt)}        &\color{mygreen}\checkmark & 0.77681478  &  0.00000000 &  6.233 &  304 \\
{\tt quadl (opt)}       &\color{mygreen}\checkmark & 0.77681478  &  0.00000000 &  3.516 &  146 \\
{\tt adaptlob (opt)}    &\color{mygreen}\checkmark & 0.77681478  &  0.00000000 &  3.515 &  164 \\
{\tt  vpaintegral}      &                         &  0.77747608  &  {\color{myfullred}0.00066129} &  1.016 &   - \\
\bottomrule
\end{tabular}
}
\endgroup
\end{table}

\begin{table}
\centering
\caption{Using optimal switching regime quadratures in Test Case 2}\label{tab:table6}
\begingroup
{\footnotesize
\begin{tabular}{lcccrr}      
\toprule
{\bf method}            & &   {\bf value}  & {\bf error}  &  {\bf time [s]} & {\bf fevals} \\
\midrule
\multicolumn{6}{l}{($\sigma=0.00001$, $o=24.447$, $f_0=0.12193$, $\omega_1=0.681$, $\omega_2=-0.914$, {\bf par = true})}  \\
\midrule
{\tt integral}           &                          & 0.00695940 &   0.00000000 &   {\color{myfullorange}5.901} &  {\color{myfullorange}312240}  \\
{\tt integral (vpa)}     &                          & 0.00695940 &   0.00000000 &  16.118 &  900 \\
{\tt integral (opt)}     &\color{mygreen}\checkmark & 0.00695940 &   0.00000000 &  13.367 &  900 \\
{\tt quad (opt)}         &                          & 0.00695930 &   0.00000010 &  12.691 &  622 \\
{\tt quadl (opt)}        &\color{mygreen}\checkmark & 0.00695940 &   0.00000000 &   8.639 &  362 \\
{\tt adaptlob (opt)}     &\color{mygreen}\checkmark & 0.00695940 &   0.00000000 &  15.986 &  686 \\
{\tt vpaintegral}        &\color{mygreen}\checkmark & 0.00695940 &   0.00000000 &  4.752  &  - \\
\midrule
\multicolumn{6}{l}{($\sigma=0.000001$, $o=28.555$, $f_0=0.12193$, $\omega_1=0.681$, $\omega_2=-0.914$, {\bf par = true})}  \\
\midrule
{\tt integral}           &                          & 0.00695451 &   {\color{myhalfred}0.00000489} &   {\color{myfullorange}10.316} &  {\color{myfullorange}561750}  \\
{\tt integral (vpa)}     &                          & 0.00695940 &   0.00000000 &  16.132 &    900 \\
{\tt integral (opt)}     &\color{mygreen}\checkmark & 0.00695940 &   0.00000000 &  13.348 &    900 \\
{\tt quad (opt)}         &                          & 0.00695930 &   0.00000010 &  12.894 &    622 \\
{\tt quadl (opt)}        &\color{mygreen}\checkmark & 0.00695940 &   0.00000000 &   8.614 &    362 \\
{\tt adaptlob (opt)}     &\color{mygreen}\checkmark & 0.00695940 &   0.00000000 &  16.104 &    686 \\
{\tt vpaintegral}        &                          & 0.00695490 &   {\color{myhalfred}0.00000449} &  22.478 &    - \\
\bottomrule
\end{tabular}
}
\endgroup
\end{table}

\subsection{Calibration to real market data}\label{ssec:cal}

The problem of calibration of the model to real market data is formulated as the nonlinear least squares optimization problem,
\[ \inf\limits_{\chi} G(\chi), \hspace{.3cm} G(\chi)=\sum^N_{i=1}w_i|C_i^\chi(T_i,K_i)-C_i^*(T_i,K_i)|^2, \]
where
$N$ denotes the number of observed option prices,
$C^*_i(T_i,K_i)$ is the market price of the call option,
$C^\chi_i$ denotes the model price of the $i$-th option computed using formula \eqref{e:price} evaluated with the vector of model parameters $\chi$, and
$w_i$ is the $i$-th weight proportional to the \emph{bid-ask spread} $\delta_i>0$.
In particular we consider the following weights \citep{MrazekPospisilSobotka16ejor}
\[
w_i = \frac{\delta_i^{-2}}{\sum\limits_{j=1}^N \delta_j^{-2}}, \quad i=1,\dots,N.
\]

In Section \ref{sec:svmodels} we mentioned that for the AFSVJD \citep{PospisilSobotka16amf,MrazekPospisilSobotka16ejor} model, the vector of model parameters is $\chi=(v_0, \kappa, \theta, \sigma, \rho, \lambda, \mu_J, \sigma_J, H)$. If not stated otherwise, the optimization will be performed with simple bounds from Table \ref{tab:bounds} and with fixed approximation parameter $\varepsilon = 10^{-6}$. Since the AFSVJD model covers also several other widely used models, we will consider also the case when $H=0.5$, i.e. the \citet{Bates96} model with 8 parameters $(v_0, \kappa, \theta, \sigma, \rho, \lambda, \mu_J, \sigma_J)$ and further with $\lambda = 0$,  i.e. the \citet{Heston93} model with 5 parameters $(v_0, \kappa, \theta, \sigma, \rho)$.

To evaluate the calibration performance, we measure the \emph{maximum and average of absolute (value of) relative error}
\[
\text{MARE}(\chi) = \max_{i=1,\dots,N} \frac{|C_i^\chi - C_i^*|} {C_i^*}, \qquad
\text{AARE}(\chi) = \frac{1}{N}\sum_{i=1}^N\frac{|C_i^\chi - C_i^*|} {C_i^*}.
\]

Following the widely accepted recommendations \citep{PospisilSobotka16amf,PospisilSobotka16dib,MrazekPospisilSobotka16ejor,MrazekPospisil17openmath} we perform the calibration as a combination of the global and local optimization techniques. Global optimization part of the calibration is especially needed when there is no suitable initial guess for the gradient based method used in the local optimizer. For the global optimizer we choose the genetic algorithm (GA)\footnote{available in MATLAB Global Optimization Toolbox, function \texttt{ga}} with the standard setting: \texttt{EliteCount} 5\% of population size, intermediate crossover (creates children by taking a random average of the parents) with \texttt{CrossoverFraction} 80\% (of the population at the next generation, not including elite children), no migration, uniform selection and Gaussian mutation. For the local optimizer we choose the standard \texttt{trust-region-reflective} Newton gradient method for nonlinear least squares (LSQ)\footnote{available in MATLAB Optimization Toolbox, function \texttt{lsqnonlin}} with stopping criteria set to function value tolerance or step tolerance \texttt{1e-9}. Calibration process is therefore performed in two steps:

\begin{description}
\item[Step 1:] Run 20 iterations (generations) of the GA optimization for utility (fitness) function $G(\chi)$, GA is run with population size 200.
\item[Step 2:] Run the LSQ optimization, the optimization is run with the initial guess obtained as a solution from the previous step.
\end{description}

In the following example, we are especially interested in comparison of calibration results when the integral in formula \eqref{e:price} is evaluated if our proposed switching algorithm is implemented or if it is not used at all. 

\begin{figure}[ht!]
\includegraphics[width=0.49\textwidth]{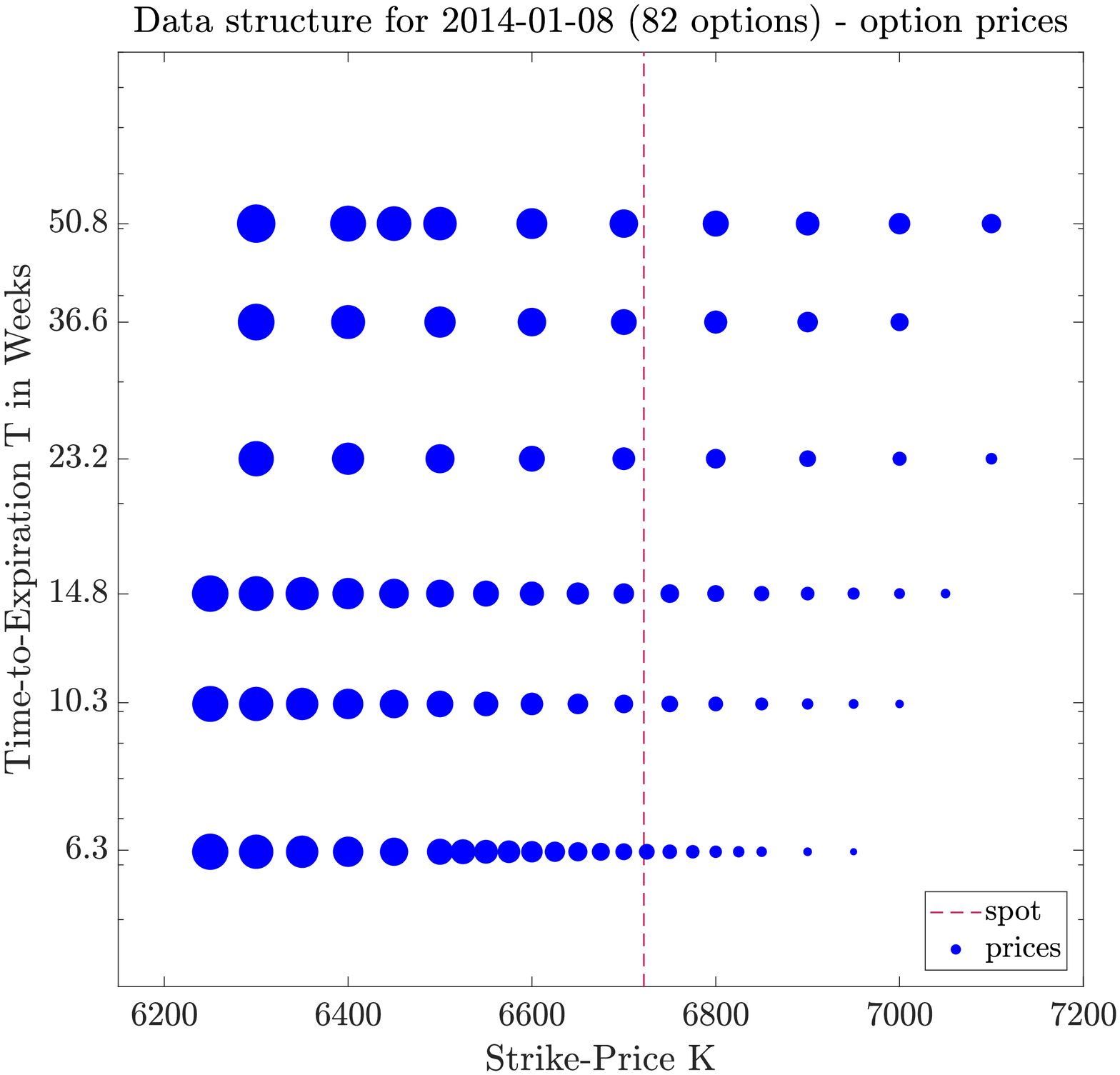}
\includegraphics[width=0.49\textwidth]{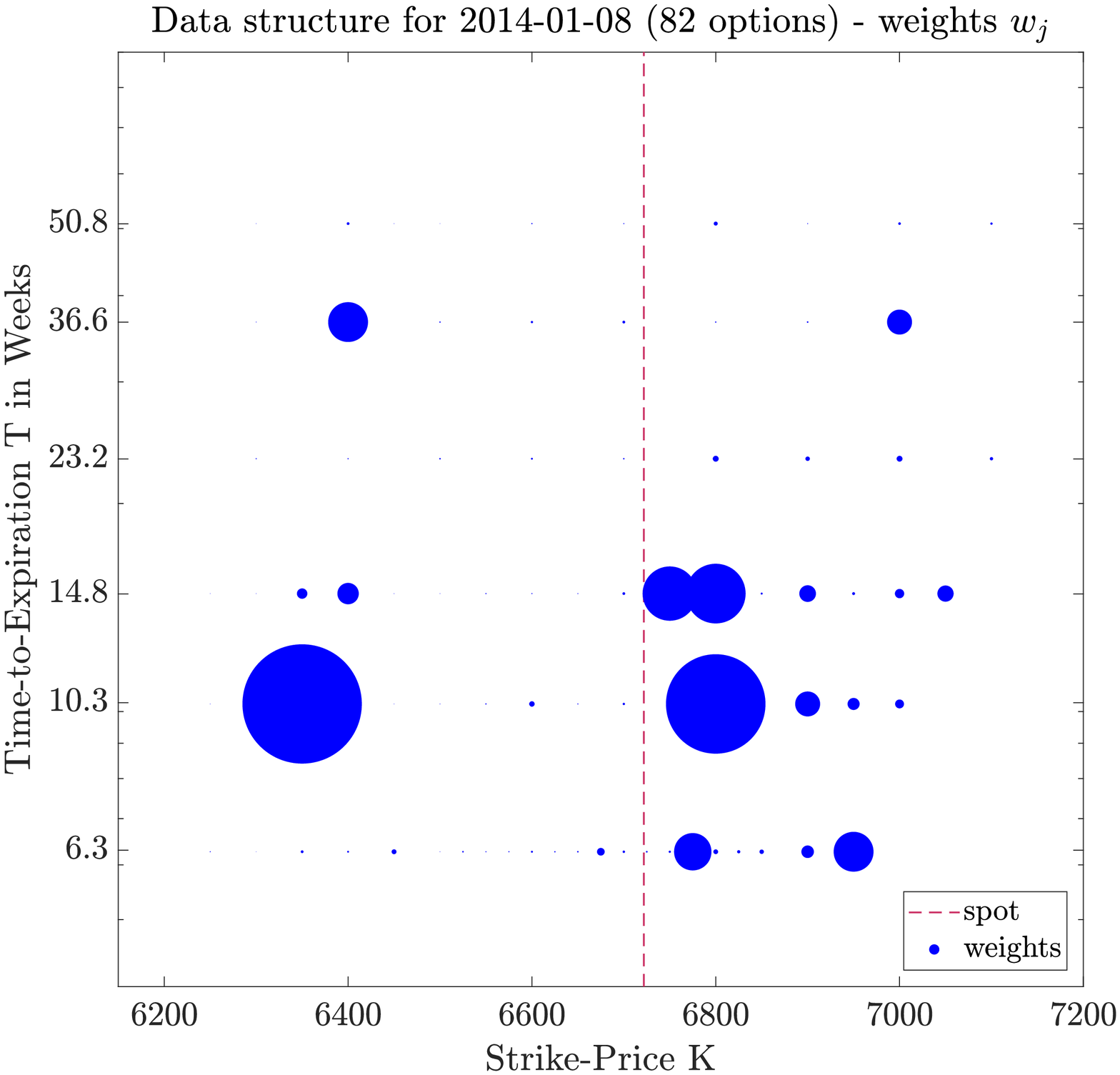}
\caption{Option data structure for data from Example \ref{ex:calib}. There are 82 options in considered data set of European call options to FTSE 100 dated 8 January 2014. Option prices are depicted on the left and considered calibration weights on the right. Note that the scaling of dot diameters differ in both pictures just for visual purposes, the diameters on the left are proportional to the option prices (USD), whereas on the right they are proportional to the unit-less weights.}\label{fig:calib-structure}
\end{figure}

\begin{example}[]\label{ex:calib}
Let us consider a data set of 82 traded European call options to the index FTSE 100 dated 8 January 2014, see \citep{PospisilSobotka16amf}. In Figure \ref{fig:calib-structure} we can see the option data structure. There are six different time to maturities ranging from $0.120548$ to $0.977528$ (in years) with strikes ranging from 6\,250 to 7\,100 with spot price 6\,721.78. In Table \ref{tab:calib-results} we can see a comparison of measured errors for all three models in two different regimes, either our fast regime switching algorithm is implemented ('ON', we use the \texttt{integral (opt)} variant from the previous section) or not used at all ('OFF'). 

As we can see, in both SVJD models (AFSVJD and Bates model) we get better calibration results if our switching algorithm is used. For Heston model we also get slightly better results, but the measured errors are indistinguishable within the tolerance \texttt{1e-9}. The table shows also the ratio how many times the evaluation of the integrand had to be switched to \texttt{vpa} from all integral calculations in each calibration task. The integral implementation is array valued, i.e. one integral calculation is performed at once for all 82 option data combinations. For convenience, we set the lower bound for parameter $\sigma$ to \texttt{1e-4} (i.e. to a value 1 bps = 0.01 \% usually used in practice). For curiosity we also include a case for AFSVJD model, where we set the lower bound for parameter $\sigma$ to \texttt{1e-2} only. As we can see, although the resulting calibrated value of parameter $\sigma$ is far from zero, avoiding the values smaller than \texttt{1e-2} during the optimization process led to a worse calibration result. On the other hand, the results 7 and 8 are almost indistinguishable within the tolerance \texttt{1e-9} since there were only 25 (in 4792) problematic situations when the integrand had to be switched to \texttt{vpa} and fortunately the final result was not affected in this case. Among all results, the 6th calibration trial gives the best results also in terms of AARE.
\end{example}

\begin{table}
\centering
\caption{Calibration results for data from Example \ref{ex:calib} for different models with switching regime algorithm (\textbf{switch}) used (ON) or not (OFF). Dimension (\textbf{dim.}) indicate the length of the vector of model parameters $\chi$ that is being calibrated. Number of switches (\textbf{\# of sw.}) indicate, how many times the evaluation of the integrand had to be switched to \texttt{vpa} in all integrand evaluations in each calibration trial. Calibrated parameter values for the AFSVJD model are listed below the table.}\label{tab:calib-results}
\begingroup
{\footnotesize
\begin{tabular}{lclrrrrr}      
\toprule
{\bf model} & {\bf dim.} & {\bf switch} & {\bf \# of sw.} & $\text{AARE}(\chi_i)$  & $\text{MARE}(\chi_i)$ &  $G(\chi_i)$ & $i$ \\
\midrule
\multirow{2}{*}{Heston} & \multirow{2}{*}{5} 
  & OFF &   0 / 4508 & 0.06512446 & 0.35292408 & 106.94072000 & 1 \\ 
 && ON  &  74 / 4442 & 0.06512469 & 0.35292146 & 106.94071999 & 2 \\ 
\midrule
\multirow{2}{*}{Bates} & \multirow{2}{*}{8} 
  & OFF &   0 / 4364 & 0.10001985 & 0.92693868 & 130.25382485 & 3 \\
 && ON  & 293 / 4526 & 0.06512521 & 0.35292381 & 106.94071992 & 4 \\
\midrule
\multirow{2}{*}{AFSVJD} & \multirow{2}{*}{9} 
  & OFF &   0 / 4812 & 0.06512390 & 0.35292162 & 106.94072035 & 5 \\
 && ON  &  88 / 6162 & 0.06065231 & 0.42843860 &  84.14079714 & 6 \\
\midrule
AFSVJD & \multirow{2}{*}{9} 
  & OFF &   0 / 4792 & 0.06512324 & 0.35291189 & 106.94071997 & 7 \\
($\sigma\ge 0.01$) &
  & ON  &  25 / 4792 & 0.06512324 & 0.35291189 & 106.94071997 & 8 \\
\bottomrule
\end{tabular}\\
\begin{align*}
\chi_5 &= (0.00644592, 6.19473908, 0.00880072, 3.95726677, 0.99999999, 0.00000000, 1.04059997, 1.61568883, 0.69621767) \\ 
\chi_6 &= (0.00001000, 0.00716856, 0.18471466, 3.71676726, 0.99995279, 3.49596838, 0.04460999, 0.00001000, 0.64278280) \\ 
\chi_7 &= (0.00644589, 6.19756427, 0.00880014, 3.23167264, 0.99999999, 0.00000246, -7.44129841, 2.73418003, 0.68154330) \\ 
\chi_8 &= (0.00644589, 6.19756427, 0.00880014, 3.23167264, 0.99999999, 0.00000246, -7.44129841, 2.73418003, 0.68154330) \\ 
\end{align*}
}\vspace*{-2em}
\endgroup
\end{table}

\section{Conclusion}\label{sec:conclusion}

In this paper we studied numerical integration in semi-closed option pricing formulas used especially in jump diffusion stochastic volatility models. When calibrating these models to real market data, a~numerical calculation of integrals has to be performed many times for different model parameters. During calibration process all integral evaluations have to be performed with high precision and low computational time requirements. Motivation for writing this paper was an observation that for some model parameters, many numerical quadrature algorithms fail to meet these requirements. We observed an~enormous increase in function evaluations (\texttt{fevals}) especially in adaptive quadratures, serious precision problems even for the simple trapezoidal rule as well as a~significant increase of computational time in all quadratures. At first we thought that the problem is in the choice of numerical quadrature. However, a more detailed numerical analysis showed that the problem is caused especially by the inaccurately evaluated integrand. We demonstrated this behaviour on a simplified integrand in conference paper \cite{DanekPospisil15tcp} and suggested the usage of variable precision arithmetic (\texttt{vpa}) in all cases when evaluating the integrand in standard \texttt{double} arithmetic is not sufficient.

The aim of this paper was to numerically analyse the integrand in the approximative fractional stochastic volatility jump diffusion model first introduced by \citet{PospisilSobotka16amf} that among others cover also the Bates and Heston models. Since the evaluation of the integrand in \texttt{vpa} is time consuming, the goal was to find a suitable fast algorithm that could tell if the integrand is \texttt{double} sufficient or if it has to be evaluated in higher precision. The main result of this paper is therefore the Algorithm \ref{alg:switch}. Experimental results then cover comparison of numerical quadratures especially for problematic (\texttt{double} insufficient) cases that has to be switched to \texttt{vpa}.

From all numerical experiments we learned several lessons. First of all, a thorough numerical analysis is a necessity in all problems that require some numerics, one should not believe any implementation of any formula if not tested thoroughly. Even a very nice looking integrand can cause serious numerical problems especially if inaccurately evaluated which can of course be hard to realize. Blindly used formulas that were not numerically analysed can lead to potential big losses as we showed in the example of wrongly priced option with difference greater than 100 dollars. We also learned that a small error of the numerical integration can be actually a coincidence and the potential problems can be detected by monitoring the number of \texttt{fevals}. And last but not least, it was not only the low value of parameter $\sigma$ that causes the problematic cases although its influence is the most remarkable.

To tell what numerical quadrature is the best for our studied integral is not an easy task, since the behaviour differs in \texttt{double} sufficient and problematic cases. In fact, when working with real data, we might find a case when every numerical quadrature fails to meet some criteria. Based on our huge number of experiments we can recommend the Gauss-Kronrod(7,15) quadrature implemented in MATLAB as a function \texttt{integral} together with the correct (\texttt{opt}) implementation of the regime switching Algorithm \ref{alg:switch}.

A~promising for the future is the idea of high-precision numerical quadrature, newly available in MATLAB function \texttt{vpaintegral} that can for example precisely calculate simple integrands such as the one in \cite{DanekPospisil15tcp}, but as we showed still has some serious problems with the integrand studied in this paper. To have a robust implementation of such a high-precision integration routine can be in fact a challenging issue, as was mentioned by \citet{Bailey05a} where several high-precision quadratures are compared on a set of test problems with integrands much simpler than the one studied here.

Last but not least, we compared the calibration results for the cases where the switching regime algorithm is used or not and showed that avoiding problematic values of some of the parameters can lead to worse calibration results.

We believe that methodology and experiments described in this paper can lead to a wider usage of high-precision numerics in financial applications and we encourage readers to perform similar numerical analysis also for other models.


\section*{Funding}

This work was partially supported by the Grantov\'{a} Agentura \v{C}esk\'{e} Republiky (GACR), grant numbers GA14-11559S \textit{Analysis of Fractional Stochastic Volatility Models and their Grid Implementation} and GA18-16680S \textit{Rough models of fractional stochastic volatility}. 

\section*{Acknowledgements}

Computational resources were provided by the CESNET LM2015042 and the CERIT Scientific Cloud LM2015085, provided under the programme ``Projects of Large Research, Development, and Innovations Infrastructure''.




\end{document}